\documentclass[11pt, reqno]{amsart}

\usepackage{amssymb}
\usepackage[T1]{fontenc}
\usepackage{pdfpages}

\evensidemargin\oddsidemargin

\begin{document}
\baselineskip=18pt
\setcounter{page}{1}
    
\newtheorem{theo}{Theorem}
\newtheorem{Lemm}{Lemma}
\newtheorem{Rem}{Remark}
\newtheorem{Coro}{Corollary}
\newtheorem{Propo}{Proposition}

\def\a{\alpha}
\def\b{\beta}
\def\A{{\bf A}} 
\def\B{{\bf B}} 
\def\C{{\bf C}} 
\def\cE{{\mathcal{E}}} 
\def\cR{{\mathcal{R}}} 
\def\cI{{\mathcal{I}}} 
\def\cB{{\mathcal{B}}}
\def\UU{{\mathcal{U}}}
\def\ca{c_{\a}}
\def\ka{\kappa_{\a}}
\def\coa{c_{\a, 0}}
\def\cua{c_{\a, u}}
\def\cL{{\mathcal{L}}} 
\def\cV{{\mathcal{V}}} 
\def\Ea{E_\a}
\def\eps{{\varepsilon}} 
\def\esp{{\mathbb{E}}} 
\def\Ga{{\Gamma}} 
\def\G{{\bf \Gamma}} 
\def\GG{{\bf G}}
\def\HH{{\bf H}}
\def\ii{{\rm i}}
\def\e{{\rm e}}
\def\L{{\bf L}}
\def\lbd{\lambda}
\def\lacc{\left\{}
\def\lcr{\left[}
\def\lpa{\left(}
\def\lva{\left|}
\def\M{{\bf M}}
\def\N{{\bf N}} 
\def\DD{{\mathbb{D}}} 
\def\NN{{\mathbb{N}}} 
\def\ZZ{{\mathbb{Z}}}
\def\pb{{\mathbb{P}}}
\def\pab{{\varphi_{a,b}}} 
\def\tpab{{{\widetilde \varphi}_{a,b}}} 
\def\rl{{\mathbb{R}}}
\def\racc{\right\}}
\def\rpa{\right)}
\def\rcr{\right]}
\def\rva{\right|}
\def\prost{{\succ_{\! st}}}
\def\T{{\bf T}}
\def\X{{\bf X}}
\def\XX{{\mathcal X}}
\def\Y{{\bf Y}}
\def\U{{\bf U}}
\def\V{{\bf V}_\a}
\def\Un{{\bf 1}}
\def\Z{{\mathbb{Z}}}
\def\A{{\bf A}}
\def\AA{{\mathcal A}}
\def\hAA{{\hat \AA}}
\def\hL{{\hat L}}
\def\hT{{\hat T}}

\def\claw{\stackrel{d}{\longrightarrow}}
\def\elaw{\stackrel{d}{=}}
\def\qed{\hfill$\square$}

\newcommand*\pFqskip{8mu}
\catcode`,\active
\newcommand*\pFq{\begingroup
        \catcode`\,\active
        \def ,{\mskip\pFqskip\relax}%
        \dopFq
}
\catcode`\,12
\def\dopFq#1#2#3#4#5{%
        {}_{#1}F_{#2}\biggl[\genfrac..{0pt}{}{#3}{#4};#5\biggr]%
        \endgroup
}

\title{On cumulative Tsallis entropies}

\author[Guillaume Dulac]{Guillaume Dulac}

\author[Thomas Simon]{Thomas Simon}

\keywords{Coherent risk measure; Cumulative entropy; Logistic distribution; Relevation process; Tsallis entropy}

\subjclass[2010]{33B15; 60E15; 62E10; 91B05; 91G70}

\begin{abstract} We investigate the cumulative Tsallis entropy, an information measure recently introduced as a cumulative version of the classical Tsallis differential entropy, which is itself a generalization of the Boltzmann-Gibbs statistics. This functional is here considered as a perturbation of the expected mean residual life via some power weight function. This point of view leads to the introduction of the dual cumulative Tsallis entropy and of two families of coherent risk measures generalizing those built on mean residual life. We characterize the finiteness of the cumulative Tsallis entropy in terms of $\cL_p-$spaces and show how they determine the underlying distribution. The range of the functional is exactly described under various constraints, with optimal bounds improving on all those previously available in the literature. Whereas the maximization of the Tsallis differential entropy gives rise to the classical $q-$Gaussian distribution which is a generalization of the Gaussian having a finite range or heavy tails, the maximization of the cumulative Tsallis entropy leads to an analogous perturbation of the Logistic distribution.  

\end{abstract}

\maketitle

\section{Introduction and notations}

Let $\DD$ be the set of real continuous random variables with finite expectation. For $X\in\DD,$ we denote by $F_X(x) = \pb[X\le x]$ the cumulative distribution function, $\bar{F}(x) = 1- F_X(x) = \pb[X > x]$ the tail distribution function and $F_X^{-1} (x) = \inf\{ y\in\rl, \; F_X(y) > x\}$ the right-continuous inverse distribution function. When there is no ambiguity we will set $F, \bar{F}$ and $F^{-1}$ respectively for $F_X, \bar{F}_X$ and $F^{-1}_X.$ In the continuous framework, one has $F(F^{-1}(x)) = x$ for every $x\in (0,1),$ the right-continuous function $F^{-1} (x)$ increases on $(0,1)$ and the positive measure $dF^{-1}(x)$ determines the law of $X,$ up to translation. Consider the so-called mean inactivity time
$${\bar \mu}(t) \; = \; \esp [t-X\,\vert \, X\le t]\; =\; \frac{1}{F(t)}\,\int_{-\infty}^t \!\! F(x)\, dx, \qquad t \in\rl,$$
with the convention that ${\bar \mu}(t) = 0$ if $F(t) = 0.$ Observe that with this convention, the mapping $t\mapsto {\bar \mu}(t)$ is a.e. differentiable, by Lebesgue's theorem. In this paper, we are interested in the following functional 
\begin{equation}
\label{Weight}
\Delta_w(X) \; =\; \esp \lcr w(F(X))\, {\bar \mu}(X)\rcr\; \in\; [0,\infty],
\end{equation}
where $w : (0,1)\to\rl^+$ is some measurable weight function. In the case of constant weight $w (x) =1$ the integration by parts formula for Riemann-Stieltjes integrals implies   
\begin{equation}
\label{Weigh}
\Delta_0(X) \; = \; - \int_\rl F(x) \log F(x)\, dx,
\end{equation}
and in the right-hand side we recognize the cumulative entropy introduced in \cite{DCL}. Here and throughout, we will set $\Delta_0(X)$ for the functional in \eqref{Weight} with weight function $w(x) = 1 = x^0.$  The cumulative entropy is a dual counterpart to the cumulative residual entropy (CRE)
\begin{equation}
\label{Dualweigh}
{\bar \Delta_0}(X) \; = \; \Delta_0(-X) \; =\; - \int_\rl {\bar F}(x) \log {\bar F}(x)\, dx\; =\; \esp \lcr {\rm mrl}(X)\rcr
\end{equation}
previously introduced in \cite{RCVW}, where 
$${\rm mrl}(t) \; = \; \esp [X-t\,\vert \, X\ge t]\; =\; \frac{1}{{\bar F}(t)}\,\int^{\infty}_t \!\! {\bar F}(x)\, dx$$
stands for the mean residual life. Observe that whereas cumulative and cumulative residual entropies can be defined for any integrable random variable, possibly taking infinite values, the identifications in \eqref{Weigh} and \eqref{Dualweigh} in terms of expected mean inactivity time resp. expected mean residual life do not hold in general if $X$ has atoms. 

The more general weight function $w(x) = (-\log x)^{n-1}/(n-1)!$ for $n\ge 1$ has also been considered in the literature, for tail distribution functions. For such weights, the same integration by part leads indeed to
$$\Delta_w(-X) \; =\; \frac{1}{n!} \int_\rl {\bar F}(x)\, (-\log {\bar F}(x))^{n} \, dx\; =\; \cE_{n} (X),\qquad X\in\DD,$$
where in the right-hand side we recognize the generalized cumulative residual entropy (GCRE) of order $n$ introduced in \cite{GN} and later investigated in \cite{DCT, Kaya, Kay, PT, TSN}, among other references. In the non-negative framework, a motivation for studying the cumulative residual entropy $\cE_1(X) = {\bar \Delta}_0(X),$ and the more general functionals $\cE_n(X),$ comes from reliability analysis and the so-called relevation transform introduced in \cite{Krav}. Specifically, if $X_0\ge 0$ denotes some failure time with reliability function ${\bar F}(t)$ and $X_1\ge 0$ denotes some second failure time with conditioned reliability function $\pb[X_1 > t \vert X_0 = x] = {\bar F} (x+t)/{\bar F} (x),$ then the "relevation process" of \cite{Krav} has lifetime $X_0+X_1,$ whose expectation is precisely $\cE_0  + \cE_1.$ More generally, the addition of further units to a given system whose failure times have the same conditioned reliability functions leads to the formula
$$\pb[T_n > t] \; =\; {\bar F}(t) \, \sum_{k=0}^{n-1} \frac{(-\log {\bar F} (t))^k}{k!}$$
where $T_n$ stands for the $n-$th failure time - see Formula (6.2) in \cite{Krav}, and hence to the identification $\esp[T_n] = \cE_0 +\cdots + \cE_{n-1}$ for the expected $n-$th failure time. Notice in passing that in \cite{Krav}, the random variable $T_n$ is considered as the instant where some wonder drug inhaled $n-$th times, loses its "full effectiveness". 

More recently in \cite{RS, CLA}, the weight function $w(x) = x^s$ for $s > -1$ was also considered for (tail) distribution functions, leading to so-called cumulative Tsallis (residual) entropies. In this framework, the integration by parts leads indeed to 
\begin{equation}
\label{Tallis}
\Delta_s(X)\; =\; \frac{1}{s} \int_\rl F(x) (1 - F(x)^s)\, dx
\end{equation}
for the cumulative entropy, and to the same formula with $F$ replaced by ${\bar F}$ for the cumulative residual entropy - see Lemma 1 in \cite{CLA}. Here and throughout, we will set $\Delta_s(X)$ for the functional in \eqref{Weight} with weight function $x^s,$ and we will make the convention $(1- x^0)/0 = -\log x$ for $x\in (0,1)$ giving \eqref{Weight} as a particular case of \eqref{Tallis}. The above terminology comes from the fact that \eqref{Tallis} can be viewed as a cumulative version of the famous Tsallis entropy
$$I_{s+1} (X) \; =\; \frac{1}{s} \lpa 1 - \int_\rl f_X^{s+1} (x)\, dx\rpa\; = \;\frac{1}{s} \int_\rl f_X(x) (1 - f_X(x)^s)\, dx$$
introduced in \cite{Tsaha} for an absolutely continuous continuous random variable $X$ with density function $f_X,$ as a generalization of the usual Shannon differential entropy corresponding to the limiting case $s = 0.$  Notice that in the case where $X$ is non-negative which is the framework of \cite{CLA}, the right-hand side of \eqref{Tallis} is well-defined and finite for every $s > -1$ and $X$ integrable, but again the identification with $\Delta_s(X)$ becomes untrue in general if $X$ has atoms, except for $s=1$ - see Remark \ref{Positano}  below. 

Let us also mention some related natural reliability model in the case $s > 0.$  Suppose that a unit in a given system has lifetime $X\ge 0$ with reliability function ${\bar F}$ and that a second unit has lifetime $Y_s\ge 0$ with conditional reliability function
$$\pb[Y_s > t \vert X= x] \; =\; {\bar F}^s (x)\lpa \frac{{\bar F}(x+t)}{{\bar F}(x)}\rpa$$
for every $t,x > 0.$ Observe that the perturbation ${\bar F}^s (x)$ implies $\pb[Y_s = 0 \vert X= x] \; =\; 1 -{\bar{F}}^s(x),$ which means that we allow for a prior failure of the second unit, with some probability increasing with both the observed value $x$ and the parameter $s.$ In the terminology of \cite{Krav}, this means that the second dosis may have lost its full effectiveness before being inhaled. Computations similar to \cite{Krav} lead then to the following expression 
\begin{equation}
\label{Rely}
\pb[X + Y_s > t] \; =\;{\bar F} (t)\lpa 1 \, +\, \frac{1- {\bar F}(t)^s}{s}\rpa
\end{equation}
for the reliability function of the lifetime of the relevation process, which implies in particular
$$\esp[Y_s] \; =\; \Delta_s(-X).$$
In the limiting case $s = 0$ without prior failure, the functional $\esp[Y_0] = \Delta_0(-X) = {\bar \Delta}_0(X)$ has been considered as a measure of dispersion or variability in Section 2 of \cite{TSN}, and Proposition 4.1 in \cite{HC} shows that the expected lifetime 
\begin{equation}
\label{Exp}
\esp[X+Y_0]\; =\; \esp[X]\, +\, {\bar \Delta}_0(X)\; =\; \esp[X + {\rm mrl} (X)]
\end{equation}
is a coherent risk measure.

The present paper goes along the previous lines of research and investigates several structural properties of the functionals $\Delta_w(X),$ highlighting those of the cumulative Tsallis entropies $\Delta_s(X).$ In Section 2, we set a sound framework for the latter functionals, extending the support from $\rl_+$ to $\rl,$ characterizing their finiteness in terms of the usual $\cL_p$ spaces, and showing how they determine the underlying distribution - see Theorem \ref{Char}. In Section 3 we introduce the dual cumulative Tsallis entropy $\nabla_s(X)$ which is given by the weight function
$$w(x)\; = \;\frac{1- (1-x)^{s+1}}{x}$$ 
and exhibits some natural duality relationship with $\Delta_s(X)$ - see Theorem \ref{Negatif}. Both cumulative and dual cumulative Tsallis entropies lead then to two natural families of coherent risk measures generalizing \eqref{Exp} and which we study thoroughly in Section 4 - see Theorem \ref{Cohe}. In Section 5 we provide several examples of distributions where the functionals $\Delta_s(X)$ and $\nabla_s(X)$ can be expressed in closed form, generalizing several computations recently made in \cite{BBL} in the case $s=0$. Finally, in Section 6, the range of $\Delta_s(X)$ acting on positive integrable distributions and on distributions with finite variance is exactly described, with sharp upper bounds improving on those recently obtained in \cite{BBL} in the case $s= 0.$ In the symmetric case, the maximizing random variables share some striking common features with the $q-$Gaussian distributions maximizing the classical Tsallis differential entropy, which were introduced in \cite{PT}. These non-standard random variables, which can be viewed as a perturbation of the Logistic, also imply a subtle inequality on a ratio of Gamma functions, which cannot seem to be simply derived by the classical  analytical arguments - see Theorem \ref{sRan} and Corollary \ref{GIneq}. Along the way and in the case $s=0$, we give a characterization of the Exponential and Logistic distributions as a maximizer in $\cL_2$ of the cumulative residual entropy resp. cumulative entropy.

\section{Cumulative Tsallis entropies} 
\label{Tsahal}

In this paragraph we consider the weight function $w(x) = x^s$ for $s > -1,$ and the corresponding functional $\Delta_s(X)$ in \eqref{Weight}. The following representation as cumulative Tsallis entropies extends Lemma 1 in \cite{CLA}, where the non-negative case was considered. In the real case, some further integrability assumption are needed on $X_- = \min(0,X)$ for $s\in (-1,0)$ in order to ensure the finiteness of $\Delta_s(X).$ The argument is standard, but we give the details for completeness. 

\begin{Propo}
\label{Affirmatif}
Let $X\in\DD.$ Then, 
$$\Delta_s(X)\; =\; \frac{1}{s} \int_\rl F(x) (1 - F(x)^s)\, dx\; \in (0,\infty]$$
for every $s \in (-1,0)\cup (0,\infty).$ For $s > 0$ one always has $\Delta_s(X) < \infty,$ whereas for $s \in (-1,0)$
$$\Delta_s(X) < \infty \;\Longleftrightarrow \; \int_{-\infty}^0 F^{1+s}(x)\, dx\; <\; \infty.$$
\end{Propo}

\proof
Suppose first $s > 0.$ Since $F(x) (1 - F(x)^s) \sim s \bar{F} (x)$ as $x\to \infty$ and $F(x) (1 - F(x)^s) \sim F(x)$ as $x\to -\infty,$ the integral on the right-hand side is finite for every $X\in\DD.$  Setting $G(x) = F(x){\bar \mu}(x),$ the integration by parts formula for Riemann-Stieltjes integrals - see e.g. Problem 6.17 in \cite{R} - implies
\begin{eqnarray*}
\frac{1}{s} \int_\rl F(x) (1 - F(x)^s)\, dx & = & \int_\rl G(x) F(x)^{s-1}\, dF(x)\, + \, \frac{\lcr G(x) (1- F(x)^s)\rcr_{\pm \infty}}{s} \\
& = & \Delta_s(X) \, + \, \lim_{x\to\infty} x {\bar F}(x) \; =\; \Delta_s(X)
\end{eqnarray*}
where the second equality comes from $G(x)\to 0$ as $x\to-\infty$ and $(1 - F(x)^s) \sim s \bar{F} (x)$ and $G(x)\sim x$ as $x\to\infty,$ and the third equality from $x {\bar F}(x)\to 0$ as $x\to\infty$ since $X$ is integrable. 

Suppose next $s\in(-1,0).$  If $\pb[X_- < c] = 0]$ for some $c\in (-\infty,0],$ we have ${\bar \mu} (x)\to 0$ as $x\to c$ and the identity
$$\Delta_s(X) \; =\; \frac{1}{s} \int_\rl F(x) (1 - F(x)^s)\, dx\; =\; \frac{1}{s} \int_c^\infty F(x) (1 - F(x)^s)\, dx\; <\; \infty$$
holds. If $X_-$ is unbounded, we have 
$$\Delta_s(X)\; = \; \frac{1}{s}\lpa \int_\rl F(x) (1 - F(x)^s)\, dx  \, - \, \lim_{x\to-\infty} G(x) F(x)^s\rpa \; \ge \; \frac{1}{s}\,\int_\rl F(x) (1 - F(x)^s)\, dx$$
and the RHS is infinite if $F^{1+s}$ is not integrable at $-\infty.$ Finally, if $F^{1+s}$ is integrable at infinity, then 
$$G(x) F^s(x)\; \le \; \int_{-\infty}^x F^{1+s} (u)\, du \;\to\; 0$$
as $x\to -\infty,$ whence
$$\Delta_s(X)\; = \; \frac{1}{s}\,\int_\rl F(x) (1 - F(x)^s)\, dx\; < \; \infty.$$
This completes the proof.

\endproof

\begin{Rem}
\label{Positano}
{\em In the case $s= 1,$ it is easy to check from the Riemann-Stieltjes integration by parts formula that the identity  
$$\Delta_1 (X)\; =\; \int_\rl G(x) \, dF(x)\, =\, \int_\rl F(x) {\bar F} (x)\, dx$$
holds true for all $X\in\cL_1$, even in the presence of atoms. Observe that in this case, we also have 
$$\Delta_1 (X)\; =\;\frac{1}{2}\, \esp [\vert X - {\tilde X}\vert]\; =\; \Delta_1 (-X)$$
where ${\tilde X}$ is an independent copy of $X.$ Observe also that for $\pb [X=0] = \pb[X=1] = 1/2,$ one has
$$\Delta_s(X) \, =\, \frac{1}{4}\qquad\mbox{and}\qquad \frac{1}{s} \int_\rl F(x) (1 - F(x)^s)\, dx\, =\, \frac{1}{2s} \, (1- 2^{-s}),$$
and that these two quantities are equal only at $s=1.$ }
\end{Rem}

The following criterion specifies the finiteness of the cumulative Tsallis entropy of index $s\in(-1,0)$ in terms of the standard $\cL_p$ spaces. Throughout, we set $\cL_p = \{ X\,\slash\,\esp[\vert X\vert^p] < \infty\}$ for every $p\ge 1.$ 

\begin{Propo}
\label{Lp}
For every $X\in\DD$ and $s\in(-1,0),$ one has the strict implications
$$X_-\in\cL_q\quad\mbox{for some $q > p$}\; \Longrightarrow\; \Delta_s(X)\, <\, \infty\;\Longrightarrow\; X_-\in\cL_p$$
with $p = 1/(1+s) \in (1,\infty).$
\end{Propo}
 
\proof

The first implication is standard since for every $q > p,$ one has
$$X_-\in\cL_q\; \Longrightarrow\; \lim_{x\to -\infty} \vert x\vert^q F(x) = 0 \;\mbox{as $x\to -\infty$} \; \Longrightarrow\; \lim_{x\to -\infty} \vert x\vert^{\frac{q}{p}} F^{1+s}(x) = 0\;\Longrightarrow\; \Delta_s(X)\, <\, \infty$$
where for the last implication we have used the last statement in Proposition \ref{Affirmatif}, which also shows that this last implication is strict. For the second implication we use the fact that $x\mapsto F^{1+s} (x)$ is a distribution function on $\rl$ for every $s \in (-1,0),$ which implies by Proposition \ref{Affirmatif}
$$ \Delta_s(X)\, <\, \infty\; \Longrightarrow\; \lim_{x\to -\infty} \vert x\vert F^{1+s}(x) = 0\;\Longrightarrow\; \vert x\vert^{p-1}  \le F^s(x) \;\mbox{for $x\in (-\infty, -M)$}$$ 
for some $M < \infty,$ whence
$$\esp[X_-^p]\; =\; \int_{-\infty}^0 \vert x\vert^p\, dF(x)\; \le\; M^p \, +\, \int_{-\infty}^0 \vert x\vert F^s(x) dF(x)\; =\; M^p\, +\, p\int_{-\infty}^0 F^{1+s} (x)\; <\; \infty$$
where the equality follows from an integration by parts. Finally, the example 
$$F(x)\; =\; \frac{1}{((1+\vert x_-\vert)\log(\e +\vert x_-\vert))^{p}}$$ shows that the second implication is strict.

\endproof 

\begin{Rem}
\label{Orli}
{\em In the case $s=0,$ one has
$$\int_{-\infty}^0 F(x)\log_+ F(x)\, dx\; <\; \infty\,\Longleftrightarrow\, \Delta_0(X)\, <\, \infty$$
for every $X\in\DD$ with the notation $\log_+ x = \sup(0,\log x).$ This yields in a similar fashion the strict implications
$$X_-\in\cL \log_+^p\! \cL\quad\mbox{for some $p > 1$}\; \Longrightarrow\; \Delta_0(X) <\; \infty\;\Longrightarrow\; X_-\in\cL\log_+\!\cL,$$
where $\cL \log_+^p\!\cL$ denotes the Orlicz space $\{ X\,\slash\,\esp[\vert X\vert\log_+^p \vert X\vert] < \infty\}$ for every $p\ge 1.$}
\end{Rem}

In the following, we will use the further notation
$$\DD_s \; =\; \{ X \in\DD, \; \Delta_s(X) < \infty\}$$
for every $s >-1.$ Observe from Propositions \ref{Affirmatif} and \ref{Lp} that this family of subsets of $\DD$ increases in $s$ with $\DD_s =\DD$ for all $s > 0.$ Recall also that all distributions in $\DD$ whose support is bounded from below are in $\DD_s$ for every $s > -1.$ The following proposition shows some useful monotonicity property for the cumulative Tsallis entropy, which holds in full generality on $\DD.$ This contrasts with the sequence $\{\cE_n, n\ge 0\},$ whose monotonicity in $n$ is connected to that of the failure rate function of $X$ in the absolutely continuous case - see equation (7) in \cite{GN}. If $X\in\DD\cap\DD_s^c$ for some $s\in (-1,0],$ we will admit infinite values for the mapping $s\mapsto\Delta_s(X)$ on $(-1,\infty)$ and set $\Delta_u(X) = \infty$ for every $u\in (-1,s].$ 

\begin{Propo}
\label{Fall}
For every $X\in\DD,$ the mapping $s\mapsto\Delta_s(X)$ decreases on $(-1,\infty)$ from $\esp[X] - \min(X)$ to $0.$ 
\end{Propo}

\proof
For every $t = F(x)\in (0,1),$ the derivative of $s \mapsto s^{-1}t(1-t^s)$ is
$$\frac{t^{s+1} (1 -t^{-s} -s \log t)}{s^2}\; < \; 0$$
on $\rl,$ and Proposition \ref{Affirmatif} implies that $s\mapsto\Delta_s(X)$ decreases on $(-1,\infty).$ By monotone convergence, the limit as $s\to\infty$ is clearly zero, whereas that as $s\to -1$ reads 
$$\esp \lcr F^{-1}(X) {\bar \mu} (X)\rcr \, =\, \int_\rl \frac{G(x)}{F^2(x)}\, dF(x) \, = \, \int_\rl {\bar F} (x) \, dx \, -\, [ {\bar F} (x) {\bar \mu} (x)]_{\pm\infty}\, =\, \int_\rl {\bar F} (x) \, dx\, + \,\lim_{x\to-\infty} {\bar \mu} (x)$$
where the second equality comes from an integration by parts and the third inequality from the asymptotic ${\bar F} (x) {\bar \mu} (x)\sim x{\bar F} (x)\to 0$ as $x\to\infty.$ If $\min(X) = -\infty$ we have
$$\esp \lcr F^{-1}(X) {\bar \mu} (X)\rcr \; \ge\; \int_\rl {\bar F} (x) \, dx\; =\; \infty,$$
whereas if $\min(X) > -\infty$ we have $\lim_{x\to-\infty} {\bar \mu} (x) = 0$ and
$$\esp \lcr F^{-1}(X) {\bar \mu} (X)\rcr \; = \; \int_\rl {\bar F} (x) \, dx\; =\; \int_0^\infty {\bar F} (x) \, dx\, -\, \int_{-\infty}^0 \!\! F (x)\, dx \, -\, \min(X) \; =\;\esp[X] \, -\, \min(X).$$

\endproof

Our next result shows the important property of cumulative Tsallis entropies that the sequence $\{\Delta_n(X), \, n\ge 1\}$ determines the law of $X\in\DD,$ up to translation. This is in sharp contrast with the cumulative entropies related to the relevation process in \cite{Krav}, where the sequence $\{\cE_n(X), \, n\ge 0\}$ may not determine the underlying distribution - see Remark \ref{Mom} below. This difference comes from the fact that the sequence $\{\Delta_n(X), \, n\ge 1\}$ defines a Hausdorff moment problem, whereas the sequence $\{\cE_n(X), \, n\ge 0\}$ defines a Stieltjes moment problem.

\begin{theo} 
\label{Char}
Let $X, Y \in\DD$ such that $\Delta_n(X) = \Delta_n(Y)$ for all $n\ge 1.$ Then, $X$ and $Y$ have the same law up to translation. 
\end{theo}

\proof Since $X\in\DD\subset\cL_1,$ the positive measure $dF^{-1}_X(x)$ on $(0,1)$ is such that
$$\int_0^1 x(1-x) dF^{-1}_X (x)\, =\, \int_\rl F_X(z) {\bar F}_X (z) \, dz\, < \, \infty$$
where the equality comes from the change of variable $z = F_X^{-1} (x),$ and the same holds for $dF^{-1}_Y(x).$ The finite positive measure $d{\tilde F}^{-1}_X(x) = x(1-x) dF^{-1}_X(x)$ clearly determines $dF^{-1}_X,$ and hence the law of $X$ up to translation. Similarly, we set $d{\tilde F}^{-1}_Y(x) = x(1-x) dF^{-1}_Y(x),$ which determines the law of $Y$ up to translation as well. By Proposition \ref{Affirmatif} and the same change of variable $z = F_X^{-1} (x),$ the equality $\Delta_n(X) = \Delta_n(Y)$ reads 
$$\int_0^1 x(1- (1- (1-x))^n)\, dF^{-1}_X(x)\; =\; \int_0^1 x(1- (1- (1-x))^n)\, dF^{-1}_Y(x)$$    
which, expanding the polynomial, implies
$$\sum_{k=1}^n (-1)^k \binom{n}{k} \int_0^1 (1-x)^{k-1} d{\tilde F}^{-1}_X(x)\; =\; \sum_{k=1}^n (-1)^k \binom{n}{k} \int_0^1 (1-x)^{k-1} d{\tilde F}^{-1}_Y(x).$$
This being true by assumption for every $n\ge 1,$ we can deduce from the triangular array that
$$\int_0^1 (1-x)^n d{\tilde F}^{-1}_X(x)\; =\; \int_0^1 (1-x)^n d{\tilde F}^{-1}_Y(x)$$
also holds true for all $n\ge 0.$ Recall that every positive finite measure $\mu$ on $(0,1)$ is determined by its integer moments since by Fubini's theorem
$$\int_0^1 e^{sx} \, d\mu (x)\; =\; \sum_{n\ge 0} \lpa \int_0^1 x^n\, d\mu(x)\rpa \frac{s^n}{n!}$$
for all $s\in\rl$ and the Laplace transform on the left-hand side characterizes $\mu.$ Putting everything together, we have shown that
$d{\tilde F}^{-1}_X(1-x) = d{\tilde F}^{-1}_Y(1-x)$ or equivalently that 
$$dF^{-1}_X(x)\, = \, dF^{-1}_Y(x)$$ 
for all $x\in (0,1),$ which completes the proof.

\endproof

\begin{Rem}
\label{Mom}{\em If $X\in\DD\cap\cL_\a$ for some $\a > 1,$ then one has $\cE_n(X) < \infty$ for every $n\ge 0.$ However, the sequence $\{\cE_n(X), \, n\ge 0\}$ may not determine the law of $X$ up to translation. Indeed, we have   
$$\cE_n (X)\; =\; \int_0^1 \frac{(-\log(1-x))^n}{n!\, x}\, d{\tilde F}^{-1}_X(1-x)\; =\; \int_0^\infty u^n\, dG_X(u)$$
with ${\tilde G}_X(u) = {\tilde F}^{-1}_X (1- e^{-u})$ and $dG_X(u) = (1-e^{-u})^{-1}d{\tilde G}_X(u)$ and by the non-uniqueness of the Stieltjes moment problem - see e.g. \cite{Lin} for a collection of criteria ensuring non-uniqueness of solutions to this problem, this implies that there exist random variables $X, Y\in\DD$ with different laws up to translation such that $\cE_n (X) = \cE_n (Y)$ for every $n\ge 0.$}
\end{Rem}

For the last result in this section we will assume that $X\in\DD$ is absolutely continuous with density function $f_X$ and we will consider its hazard rate function
$$\lbd_X(t) \; =\;\frac{f_X(t)}{{\bar F}_X(t)}$$ 
with the convention that $\lbd_X(t) = 0$ for $t\ge \max X.$ Again we will set $f, \lbd$ for $f_X, \lbd_X$ respectively when there is no ambiguity. Following Chapter 1 in \cite{SS}, we will say that $X$ is DFR (decreasing failure rate) if $\lbd$ is non-increasing on Supp $X,$ in other words if $\log {\bar F}$ is convex. Observe that if $X$ is DFR, then Supp $X$ must be an interval and that $f$ must be positive on the interior of this interval. Assuming furthermore that $X$ is non-negative and setting $X_s = X + Y$ for the lifetime of the relevation process in the introduction with $s \ge 0,$ it follows from \eqref{Rely} that $X_s$ has density function 
$$f_s (t) \; =\; \frac{(1+s)\, f(t)\, (1-{\bar F}(t)^s)}{s}\cdot$$
Setting, here and throughout, ${\bar \DD}_s =\{ X\in\DD\,\slash\, -\!X \in\DD_s\}$ and ${\bar \Delta}_s(X) = \Delta_s(-X)$ for $X\in{\bar \DD}_s,$ we can deduce similarly as in \cite{GN} that
$${\bar \Delta}_s(X) \; =\; \frac{1}{s} \int_\rl {\bar F}(x) (1 - {\bar F}(x)^s)\, dx\; =\; \frac{1}{s} \int_\rl (1 - {\bar F}(x)^s)\, \frac{f(x)}{\lbd (x)}  \,dx\; =\; \esp\lcr \frac{1}{\lbd(X_s)}\rcr.$$
A consequence of this representation is the following ordering result echoing Theorem 1 in \cite{GN}. Recall from Chapter 1.B.1 in \cite{SS} that for two absolutely continous random variables $X$ and $Y,$ the hazard rate ordering $$X\,\preceq_{hr}\, Y$$ 
means $\lbd_X(x)\ge \lbd_Y(x)$ for all $x\in\rl,$ in other words the mapping $x\mapsto {\bar F}_X(x)/{\bar F}_Y(x)$ is non-increasing. 

\begin{Propo}
\label{HR}
For every $s > -1,$ if $X, Y\in\DD_s$ with $X\preceq_{hr} Y$ and either $X$ or $Y$ is {\rm DFR}, then one has
$${\bar \Delta}_s(X) \; \le\;{\bar \Delta}_s(Y).$$
\end{Propo}

\proof 
 It follows from \eqref{Rely}, the increasing character of the function
 $$t\;\mapsto\; t\, +\, \frac{t(1-t^s)}{s}$$
 on $(0,1)$ with derivative $(s+1)(1-t^s)/s > 0,$ and Theorem 1.B.1 in \cite{SS} that 
$$X\,\preceq_{hr}\, Y\;\Longrightarrow\; X\,\preceq_{st}\, Y  \;\Longrightarrow\; X_s\,\preceq_{st}\, Y_s$$
where, here and throughout, the stochastic ordering $X\,\preceq_{st}\, Y$ means ${\bar F}_X(x)\ge {\bar F}_Y(x)$ for all $x\in\rl.$ The proof goes then along the same lines as in Theorem 1 of \cite{GN}.
   
\endproof

\section{Dual cumulative Tsallis entropies} 
\label{Duelle}

In this paragraph we consider the weight function
$$w_s(x) = \frac{1 -(1-x)^{s+1}}{x}$$
for $s > -1,$ and we call $\nabla_{\! s}(X)$ the corresponding functional in \eqref{Weight}. At first sight, this weight function looks more complicated than the Tsallis weight $x^s$ investigated in the previous section. Observe however that since $w_s(x) \to (s+1) > 0$ as $x\to 0$ the finiteness of $\nabla_{\! s}(X)$ for $X\in \DD,$ which is read off on the sole behaviour of $F(x)$ at $-\infty$ since $w_s(x) \to 1$ as $x\to 1,$ amounts to that of the case with constant weight. In other words, one has 
$$\nabla_{\! s}(X)\, <\,\infty\;\Longleftrightarrow\; X\,\in\,\DD_0$$
for every $X\in\DD$ and $s > -1.$ Also, the increasing character of $s\mapsto w_s(x)$ implies as in Proposition \ref{Fall} that for every $X\in\DD_0,$ the mapping $s\mapsto \nabla_{\! s}(X)$ increases on $(-1,\infty)$ from $0$ to $\esp[X] - \min(X).$

\medskip

The following result exhibits an interesting duality relationship between $\Delta_s(X)$ and $\nabla_{\! s}(X).$ Throughout, we will use the Pochhammer notation $(x)_0 = 1$ and $(x)_n = x(x+1)\ldots (x+n-1), \,n \ge 1,$ for the ascending factorial of a real number $x.$

\begin{theo}
\label{Negatif} {\em (a)} For every $s > -1$ and $X\in\DD_0,$ one has the convergent series representation
$$\nabla_{\! s}(X) \; =\; (1+s) \sum_{n\ge 0} \frac{(-s)_n}{(n+1)!}\,\Delta_{n}(X).$$

{\em (b)} For every $s > -1$ and $X\in\DD_s,$ one has the convergent series representation
$$\Delta_s(X)\; =\; (1+s) \sum_{n\ge 0} \frac{(-s)_n}{(n+1)!}\,\nabla_{\! n}(X).$$
\end{theo}

\proof Expanding by the generalized binomial theorem $1 - (1-F(x))^{s+1}$ with $F(x)\in (0,1),$ we first compute 
$$w_s(F(x))\; =\; - \sum_{n\ge 1} \frac{(-s-1)_n}{n!}\, F^{n-1} (x) \; =\; (1+s) \sum_{n\ge 0}\frac{(-s)_n}{(n+1)!}\, F^n(x).$$
The terms of the above series have constant sign for $n$ large enough, and  Tonelli's theorem  combined with the definition $\nabla_{\! s}(X) = \esp[w_s(X){\bar \mu}(X)]$ implies
\begin{equation}
\label{plusmin}
\nabla_{\! s}(X) \; = \;  (1+s) \sum_{n\ge 0} \frac{(-s)_n}{(n+1)!}\,\Delta_{n}(X)
\end{equation}
for every $X\in \DD_0$ and $s > -1$ as required for (a). The proof for (b) is analogous, using 
 \begin{eqnarray*}
x^s \; = \; \frac{(1- (1-x))^{s+1}}{x} \; = \; \frac{1}{x}\, \sum_{n\ge 0}\frac{(-s-1)_n\,(1-x)^n}{n!}
& = & \frac{1}{x}\, \sum_{n\ge 0}\frac{(-s-1)_n\,((1-x)^n -1)}{n!}\\ & = & (s+1) \sum_{n\ge 0}\frac{(-s)_n}{(n+1)!}\lpa\frac{1- (1-x)^{n +1}}{x}\rpa
\end{eqnarray*}
for every $x\in (0,1),$ where in the third equality we have used
$$\sum_{n\ge 0}\frac{(-s-1)_n}{n!}\; = \; 0$$
since $-s-1 < 0.$ 
\endproof

\begin{Rem}
\label{Senans}
{\em (a) Part (a) of the previous result gives an interpretation of $\nabla_{\! s}(X)$  as the expected lifetime of some randomized relevation process, for $s\in(-1,0)$. Indeed, setting $A_u$ for an arcsine random variable with parameter $u = -s \in (0,1)$ and probability density function
$$\frac{\sin \pi u}{\pi}\; x^{-u} (1-x)^{1-u}\, \Un_{(0,1)} (x),$$
and $G_p$ for an independent geometric random variable with parameter $p\in (0,1)$ and probability mass function $\pb[G_p = n] = (1-p)p^n$ on $\NN,$ a computation on generating functions shows that the randomized random variable $N_s = G_{A_s}$ has probability mass function
$$\pb[N_s = n] \; = \; (1+s)\, \frac{(-s)_n}{(n+1)!},\quad n\ge 0.$$
In the literature, the law of the random variable $N_s$ is sometimes called a beta negative binomial distribution ${\rm BNB} (r,\a,\beta)$ with parameters $r=1, \a = 1+s$ and $\beta = -s.$ If we now set $X_n = X+Y$ for the lifetime of the relevation process discussed in the introduction with integer parameter $n,$ and if we consider an independent random variable $N_s,$ we see that 
$$\nabla_{\! s} (X)\; = \; \esp\lcr X^{}_{N_s}\rcr\!.$$
In the case $s > 0$ however, a direct connection between $\nabla_{\! s}(X)$ and relevation processes is less clear because of the negative weights. Notice that in the case $s\in (-1,0),$ there is no clear connection between $\Delta_s(X)$ and relevation processes  either, since the weight ${\bar F} (x)^s$ becomes greater than one. An interesting open problem is whether there exists other randomizations depending on $s$ leading to duality formulas as in Theorem \ref{Negatif}.

\medskip

(b) If we take $s \in (-1,0)$ in Theorem \ref{Negatif} (b), the well-known asymptotic
$$\frac{(-s)_n}{(n+1)!}\, \sim\, \frac{1}{\Gamma(-s) n^{2+s}}\qquad \mbox{as $n\to\infty$}$$
implies 
$$\sum_{n\ge 0} \frac{\nabla_n(X)}{n^{2+s}}\; <\; \infty$$
for every $X\in\DD_s.$ In the case when $X_-$ is unbounded, this gives some information on the speed of convergence of $\nabla_n(X)$ towards $\infty.$ It would be interesting to study this speed more precisely, as well as that of $\Delta_n(X)$ towards zero. A first observation in this respect is that
$$n \Delta_n(X)\; \longrightarrow \; \int_\rl F(x)\, dx\; =\; \max(X) \, -\, \esp[X]$$
as $n\to\infty$ and that the RHS is infinite if $X_+$ is unbounded.

\medskip
 
(c) If we take $s=k$ an integer in Theorem \ref{Negatif} (a) and (b), we obtain 
$$\nabla_k(X) \; = \;  \sum_{n=0}^k \binom{k+1}{n+1}\, (-1)^n\, \Delta_{n}(X) \qquad\mbox{and}\qquad\Delta_{k}(X) \; = \;  \sum_{n=0}^k \binom{k+1}{n+1}\, (-1)^n\,\nabla_n(X)$$
for every $k\ge 0,$ which corresponds to the easily established involutory character of a certain family of alternate binomial matrices.}
\end{Rem}

Putting together Theorems \ref{Char} and \ref{Negatif} (b) gives the following  characterization property of the dual cumulative Tsallis entropy.

\begin{Coro} For every $X\in\DD_0,$ the sequence $\{\nabla_n(X),\, n\ge 0\}$ determines the law of $X$ up to translation.
\end{Coro}

We end this section with the introduction of two new skewness parameters defined in terms of the operators $\Delta_s$ and $\nabla_s,$ and which we believe to be relevant. Recalling the notations ${\bar \Delta}_s (X) = \Delta_s(-X)$ and ${\bar \DD}_s = \{X\in\DD\,\slash -\!X \in\DD_s\}$ and setting, here and throughout, ${\bar \nabla}_s (X) = \nabla_s(-X),$ it follows from the above discussion that the two mappings
$$\Diamond_X : s\, \mapsto\, \frac{\nabla_s(X)}{{\bar \Delta}_s(X)}\qquad\quad\mbox{and}\qquad\quad {\bar \Diamond}_X : s\, \mapsto\, \frac{{\bar \nabla}_s(X)}{\Delta_s(X)}$$
are well-defined and right-continuous increasing from 0 to $\infty$ on $(-1,\infty)$ for every $X\in\DD_0\cap {\bar \DD}_0.$ Moreover, they are continuous except for possibly one jump on $(-1,0]$ from the value $0$. Defining 
$$\varrho(X) = \inf\{ s > - 1\; \slash\, \Diamond_X(s) > 1\}\qquad\quad \mbox{and}\qquad\quad {\bar \varrho}(X) = \inf\{ s > - 1\; \slash\, {\bar \Diamond}_X(s) > 1\},$$
we see that $\varrho(X) = {\bar \varrho} (X) = 0$ if and only if  $\Delta_0(X) = {\bar \Delta}_0(X),$ which happens e.g. if $X$ is symmetric around its mean. The following specifies the range of $\varrho(X)$ and ${\bar \varrho}(X)$ when $\Delta_0(X) \neq {\bar \Delta}_0(X).$

\begin{Propo} Let $X\in \DD_0\cap {\bar \DD}_0.$ If $\Delta_0(X) > {\bar \Delta}_0(X),$ then one has $-1 < \varrho(X) \le 0 < {\bar \varrho} (X) < 1.$ If $\Delta_0(X) < {\bar \Delta}_0(X),$ then one has $-1 < {\bar \varrho}(X) \le 0 < \varrho (X) < 1.$
\end{Propo}

\proof It suffices to consider the case $\Delta_0(X) > {\bar \Delta}_0(X).$ Since $\nabla_0(X) = \Delta_0(X)$ we have $\Diamond_X(0) > 1$ and hence $\varrho(X) \le 0.$ Since ${\bar \nabla}_0(X) = {\bar \Delta}_0(X)$ we have ${\bar \Diamond}_X(0) < 1$ and hence ${\bar \varrho}(X) > 0$ by continuity of $s\mapsto {\bar \Diamond}_X(s)$ on $[0,\infty).$ Finally, since ${\bar \nabla}_1(X) > {\bar \Delta}_1(X) = \Delta_1(X)$ by Remark \ref{Positano}, we have ${\bar \Diamond}_X(1) > 1$ and hence ${\bar \varrho}(X) < 1$ by the same continuity of $s\mapsto {\bar \Diamond}_X(s)$ on $[0,\infty).$

\endproof

It is interesting to mention - see Remark \ref{RAB} below - that the mappings $X\mapsto\varrho(X)$ and $X\mapsto{\bar \varrho}(X)$ may take any value in $(-1,1),$ which is a nice property for a statistical parameter. The extremal values - 1 and 1 are approached by distributions close to Dirac measures on the left resp. on the right. In general, the parameters $\varrho(X)$ and ${\bar \varrho}(X)$ are related to the dispersion and the asymmetry of $X.$ Their further properties will be investigated elsewhere. 
  
\section{Two families of coherent risk measures} 
\label{Coco}

In this paragraph we consider the functionals 
$$\nu_s^\Delta (X) \; =\; \esp[X]\; +\; {\bar \Delta_s(X)}\qquad\quad \mbox{and}\qquad\quad \nu_s^\nabla (X) \; =\; \esp[X]\; +\; {\bar \nabla_s(X)}$$
for $s > -1.$  The functional $\nu_s^\Delta$ takes finite values on ${\bar \DD}_s,$ whereas the functional $\nu_s^\Delta$ takes finite values on ${\bar \DD}_0,$ for every $s >-1.$ Observe that for $s = 0$ we also have 
$$\nu_0^\Delta (X) \; =\; \nu_0^\nabla (X) \; =\; \esp[X +{\rm mrl} (X)]$$
with the notation ${\rm mrl} (X)]$ of the introduction for the mean residual life of $X$. The following alternative representation shows how $\nu_s^\Delta (X)$ and $\nu_s^\nabla (X)$ can be viewed as perturbations of the functional $\esp[X +{\rm mrl} (X)]$, and also gives a new aspect of the duality relationship in Theorem \ref{Negatif}.

\begin{Propo}
\label{DN}
For all $s > -1,$ one has
$$\nu_s^\Delta (X) \; =\; (s+1)\, \esp [{\bar F}(X)^s ( X +{\rm mrl} (X))]\qquad\mbox{and}\qquad \nu_s^\nabla (X) \; =\; (s+1)\, \esp\lcr F(X)^s_{} ( X +{\rm mrl} (X))\rcr$$
for $X\in {\bar \DD}_s$ resp. $X\in {\bar \DD_0}.$
\end{Propo}

\proof

Setting $h(x) = x + {\rm mrl} (x),$ we compute $dh(x) = {\bar F} (x)^{-1} {\rm mrl} (x) dF(x).$ The Riemann-Stieltjes integration by part formula yields
\begin{eqnarray*}
(s+1)\, \esp [{\bar F}(X)^s ( X +{\rm mrl} (X))] & = & (s+1) \,\int_\rl {\bar F}(x)^s h(x)\, dF(x)\\
& = & \int_\rl {\bar F}(x)^{s+1}\, dh(x)\, -\, \lcr h(x) {\bar F}(x)^{s+1}\rcr_{\pm\infty}\\
& =& \int_\rl {\bar F}_X(x)^s\,{\rm mrl}_X (x) \, dF_X(x)\, +\, \esp[X]\\
& =& \int_\rl F_{-X}(-x)^s \,{\bar \mu}_{-X} (-x) \, dF_{-X}(-x)\, +\, \esp[X]\; = \; \nu_s^\Delta (X)
\end{eqnarray*}
where in the fourth equality we have used ${\rm mrl}_X (x) = - {\bar \mu}_{-X} (-x)$ for all $x\in\rl$ and the definition of $\Delta_s(-X),$ whereas the third equality comes from the further integration by parts 
$$h(x) {\bar F}(x) \; =\; x{\bar F}(x)\, +\, \int_x^\infty {\bar F}(t)\, dt \; =\; \int_x^\infty t\, dF(t)$$
which gives $h(x) {\bar F}(x)^{s+1} \to \esp [X]$ as $x\to -\infty$ for all $s > -1,$   $h(x) {\bar F}(x)^{s+1} \to 0$ as $x\to \infty$ for all $s\ge 0,$ and
$$h(x) {\bar F}(x)^{s+1}\; \le \; x{\bar F}(x)^{s+1}\, +\, \int_x^\infty {\bar F}(t)^{s+1}\, dt \; \to\; 0$$
as $x\to \infty$ for all $s\in (-1,0)$ and $X\in {\bar \DD}_s,$ by Proposition \ref{Affirmatif}. Similarly, we have
\begin{eqnarray*}
(s+1)\, \esp [F(X)^s ( X +{\rm mrl} (X))] & = & (s+1) \,\int_\rl F(x)^s h(x)\, dF(x)\\
& = & \int_\rl (1- F(x)^{s+1})\, dh(x)\, +\, \lcr h(x) (F(x)^{s+1} -1)\rcr_{\pm\infty}\\
& =& \int_\rl w_s({\bar F}_X(x))\,{\rm mrl}_X (x) \, dF_X(x)\, +\, \esp[X]\; = \; \nu_s^\nabla (X)
\end{eqnarray*}
where the identification  $\lcr h(x) (F(x)^{s+1} -1)\rcr_{\pm\infty} = \esp[X]$ is derived as above and is here valid for $X\in {\bar \DD}_0.$ We omit details.
\endproof

We now establish some interesting connections between the two above functionals and risk theory.   Following e.g. \cite{Ch} p. 5, we recall that a functional $\cR$ acting on a given space of real random variables is a coherent risk measure if it satisfies the following properties:

\begin{itemize}
\item Law-invariance: $X\elaw Y \,\Rightarrow \,\cR(X) = \cR (Y).$
\item Monotonicity: $X\preceq_{st} Y \Rightarrow \cR(X) \le \cR (Y).$
\item Translation invariance: $\cR(X +b) = \cR(X) + b $ for all $b\in\rl.$
\item Homogeneity: $\cR(aX) = a\cR(X)$ for all $a > 0$
\item Subadditivity: $\cR(X+Y) \le \cR(X) +\cR(Y).$
\end{itemize}

Above, the subadditivity property is meant without independence assumptions on $X$ and $Y$ and the stochastic order $\preceq_{st}$ is defined as usual by $X\preceq_{st} Y \Leftrightarrow {\bar F}_X\le {\bar F}_Y$ on $\rl.$ As mentioned in the introduction, it was shown in Proposition 4.1 in \cite{HC} that the functional $\esp[X +{\rm mrl} (X)]$ is a coherent risk measure. The following result extends this property to all $s >-1.$

\begin{theo}
\label{Cohe}
The functionals $\nu^\Delta_s$ and $\nu^\nabla_s$ are coherent risk measures on ${\bar \DD}_s$ resp. on ${\bar \DD_0}.$  
\end{theo}

\proof

Law-Invariance is clear, and it follows from the definitions that ${\bar \Delta}_s(aX +b) = a{\bar \Delta}_s(X)$ and ${\bar \nabla}_s(aX +b) = a{\bar \nabla}_s(X)$ for all $a > 0, b,s\in\rl$ and $X\in{\bar \DD}_s$ resp. $X\in {\bar \DD_0},$ which yields homogeneity and translation invariance. To obtain the monotonicity and subadditivity properties, we will prove that the functionals $\nu_s^\Delta$ and $\nu_s^\nabla$ are so-called Wang distortion measures. This property is easy to see for $\nu_s^\Delta$ since
\begin{eqnarray*}
\nu^\Delta_s(X)& = &  \int^{\infty}_0 {\bar F}(x)\, dx\, -\, \int_{-\infty}^0 (1- {\bar F}(x))\, dx\, +\, \frac{1}{s} \,\int_\rl {\bar F}(x) (1- {\bar F}^s(x))\, dx\\
& = & \int^{\infty}_0 h_s({\bar F}(x))\, dx\, -\, \int_{-\infty}^0 (1- h_s({\bar F}(x)))\, dx
\end{eqnarray*}
with 
$$h_s(t)\; =\; t \, +\, \frac{t(1-t^s)}{s}$$
for all $t\in (0,1).$ Computing $h_s'(t) = (s +1)(1-t^s)/s > 0$ on $(0,1)$ as in Proposition \ref{HR} shows that $\nu_s$ is monotone. Computing further $h_s''(t) = -(s+1) t^{s-1} < 0$ on $(0,1)$ shows that $h_s$ is concave, which implies by Theorem 10 in \cite{WD} that $\nu_s^\Delta$ is subadditive.

The increasing concave distortion property for $\nu_s^\nabla$ requires some more effort. By Theorem \ref{Negatif} (a), we have 
$$\nu^\nabla_s(X)\; = \; \esp[X] \, +\, {\bar \Delta}_0(X) \, +\,(1+s) \sum_{n\ge 1} \frac{(-s)_n}{(n+1)!}\,{\bar \Delta}_n(X),$$
which leads by Tonelli's theorem to the decomposition 
$$\nu_s^\nabla(X)\; =\; \int^{\infty}_0 k_s({\bar F}(x))\, dx\; -\; \int_{-\infty}^0 (1- k_s({\bar F}(x)))\, dx$$
with $k_s(t) \, = \,  t \lpa 1  \, +\, (1+s) \lpa K_s(1) \, -\, K_s(t)\, - \, \log t\rpa\rpa$
and the further notation 
$$K_s(t)\; =\; \sum_{n\ge 1} \frac{(-s)_n\, t^n}{n(n+1)!}\cdot$$
We finally compute
$$k_s'(1)\; =\; -(s+1) K_s'(1)\, -\, s \;  = \;  -(s+1)\sum_{n\ge 1} \frac{(-s)_n}{(n+1)!}\, -\, s \; = \; 0$$
and
\begin{eqnarray*}
k_s''(t)\, = \, -(s+1)\lpa t^{-1}  + 2K_s'(t) + tK_s''(t)\rpa &  = & -(s+1)\, \sum_{n\ge 0}\frac{(-s)_n t^{n-1}}{n!} \\
& = & -(s+1)\, t^{-1} (1-t)^s\; <\; 0,
\end{eqnarray*}
which also implies $k_s'(t) > 0$ on $(0,1)$.

\endproof

\begin{Rem}{\em (a) The homogeneity and subadditivity properties of $\nu_s^\Delta$ and the linearity of the expectation imply the inequality 
$${\bar \Delta}_s(\lbd X + (1-\lbd) Y)\; \le\; \lbd {\bar \Delta}_s(X) \, +\, (1-\lbd) {\bar \Delta}_s(Y)$$
for all $s\in \rl, \lbd\in [0,1]$ and $X,Y$ integrable, which is known as a convexity property of a risk measure. In the case $s = 0$ of the cumulative residual entropy, the latter inequality was claimed in Proposition 2.7 of \cite{TP} as a consequence of the representation ${\bar \Delta}_0(X) = \esp[T(X)]$ where $T$ is a certain convex function, which depends however on the distribution function of $X$ so that this argument is incomplete. Observe that neither ${\bar \Delta}_s(X)$ nor ${\bar \nabla}_s(X)$ itself is a coherent risk measure, and that the addition of $\esp[X]$ is crucial to make both functionals monotonous.

\medskip

(b) By Proposition \ref{Fall} the mapping $s\mapsto \nu^\Delta_s(X)$ decreases on $(-1,\infty)$ from $\max X$ to $\esp[X],$ whereas by the discussion at the beginning of Section \ref{Duelle} the mapping $s\mapsto \nu^\nabla_s(X)$ increases on $(-1,\infty)$ from $\esp[X]$ to $\max X.$ Both functionals $\max X$ and $\esp[X]$ are well-known, basic examples of coherent risk measures.

\medskip

(c) The above proof shows that $X\preceq_{st} Y\Rightarrow \nu^\Delta_s(X) \le \nu^\Delta_s(Y)$ and $\nu^\nabla_s(X) \le \nu^\nabla_s(Y)$ for all $s > -1.$ It is worth mentioning that the converse is not true in general. For example, if $X_{a,L}$ is uniformly distributed on $(a,a+L),$ then some computations using Propositions \ref{Affirmatif} and \ref{Negatif} - see also Paragraphs \ref{Beta1} or \ref{Beta2} with $\beta = 1$ below - give
$$\nu_s^\Delta(X_{a,L}) \; =\; a\, +\, \frac{L(s+3)}{2(s+2)}\qquad\mbox{and}\qquad \nu_s^\nabla(X_{a,L}) \; =\; a\, +\, \frac{L(2s+3)}{2(s+2)}$$
for all $s > -1.$ Hence, if $b > a$ and $L = M + 2(b-a)$ then $\nu_s^\Delta(X_{b,M}) < \nu_s^\Delta(X_{a,L})$ and $\nu_s^\nabla(X_{b,M}) < \nu_s^\nabla(X_{a,L})$ for every $s > -1$ but clearly there is no stochastic ordering between $X_{a,L}$ and $X_{b,M}$ since $a < b < b + M < a+L.$ 

\medskip

(d) For $s \ge 0,$ consider the functional
$${\tilde \nu}_s(X)\; = \; \esp[X] \; +\; \cE_s(X)$$
which is built on the generalized cumulative residual entropy 
$$\cE_s(X)\; =\; \frac{1}{\Ga(s+1)}\int_\rl {\bar F} (x) \lpa -\log {\bar F} (x)\rpa^s \, dx$$
introduced in \cite{GN}. One has ${\tilde \nu}_0(X) = 2\esp[X] - \min X$ and ${\tilde \nu}_1(X) = \nu_1(X),$ which are both coherent risk measures. However, for $s\not\in\{0,1\}$ one has
$${\tilde \nu}_s(X)\; =\; \int^{\infty}_0 {\tilde h}_s({\bar F}(x))\, dx\, -\, \int_{-\infty}^0 (1- {\tilde h}_s({\bar F}(x)))\, dx$$
with 
$${\tilde h}_s(t) \; =\; t\, +\, \frac{t \,(-\log t)^s}{\Ga(s+1)}$$
and one can check that ${\tilde h}_s$ is not non-decreasing on $(0,1)$ for $s \in (0,1)$ and not concave on $(0,1)$ for $s > 1$ - see also Remark 4.2 in \cite{HC} for the case when $s$ is an integer. By the main result in \cite{Sc} this implies that ${\tilde \nu}_s$ is not monotonous for $s\in (0,1),$ and by Proposition 3 in \cite{Sc} that it is not subadditive for $s > 1.$ In particular, the functional ${\tilde \nu}_s(X)$ is a coherent risk measure for $s\in\{0,1\}$ only, in contrast to Theorem \ref{Cohe}. On the other hand, the function
$${\tilde H}_n(t) \; =\; \sum_{k=0}^n \frac{t \,(-\log t)^n}{n!}$$
has derivative $(-\log t)^n/n!$ which is positive and decreasing on $(0,1),$ so that the functional $\esp[T_n],$ with the notation of the introduction, is a coherent risk measure for all $n\ge 1.$}
\end{Rem}

\section{Some explicit examples}

In this paragraph we display some random variables whose cumulative Tsallis entropies and dual cumulative Tsallis entropies can be computed in closed form, mostly in terms of the Gamma function $\Gamma$ and the Digamma function $\psi$, which we recall to be defined as
$$\psi(z)\; =\; -\gamma + \sum_{n\ge 0} \lpa \frac{1}{n+1} - \frac{1}{n+z}\rpa\; =\; - \gamma + \int_0^1 \lpa\frac{1- t^{z-1}}{1 - t}\rpa dt, \qquad z > 0,$$
where $\gamma$ is the Euler-Mascheroni constant. The list is not exhaustive, and for the sake of concision we will also not give the full details behind the computations. Some formulas will be used in the next section when investigating the range of $\Delta_s.$ We will mostly consider explicit transformations of $\U,$ the uniform random variable on $(0,1)$ and $\L,$ the standard exponential random variable with density $e^{-x}$ on $(0,\infty).$ We recall that for all $s >-1, a > 0$ and $b\in\rl,$ one has
\begin{equation}
\label{Aff}
\Delta_s (aX + b)\; =\; a \Delta_s(X) \qquad\mbox{and}\qquad \nabla_{\!s}(aX+b)\; =\; a \nabla_{\!s}(X)
\end{equation}
for $X\in\DD_s$ resp. $X\in\DD_0.$ On the other hand, there is no such simple relationship for $a <0$ in general, except for $\Delta_1$ with $\Delta_1 (aX + b) = \Delta_1(a X) = -a \Delta_1(X)$ - see Remark \ref{Positano}. 

\subsection{$X = \U^{1/\beta}, \, \beta > 0$} 

\label{Beta1}

The density function is $f(x) = \beta x^{\beta -1} \Un_{(0,1)} (x).$ For all $s> -1,$ we have
$$\Delta_s = \frac{\beta}{(\beta + 1)(\beta (1+s) +1)}\qquad\mbox{and}\qquad \nabla_{\! s} = \frac{\beta}{\beta + 1}\lpa 1 - \frac{\Gamma(1/\beta +1)\Gamma(s+2)}{\Gamma(1/\beta + s +2)}\rpa.$$
Observe, with obvious notations, that 
$$\nabla_{-1}\, =\, \Delta_\infty \,=\, 0 \qquad \mbox{and}\qquad \Delta_{-1} \, =\, \nabla_\infty \,=\, \frac{\b}{\b+1}\, = \,\esp[X]$$ 
as expected. In the case $s = 0$ of the cumulative entropy, we get 
$$\Delta_0 \; = \; \nabla_{\! 0} \; =\; \frac{\beta}{(\beta +1)^2}\cdot$$
In particular, the case $\beta = 1$ yields $\Delta_0 = 1/4,$ which was recently evaluated in Example 2 of \cite{BBL} by other methods.

\subsection{$X = 1 - \U^{1/\beta}, \, \beta > 0$} 

\label{Beta2}

The density function is $f(x) = \beta (1-x)^{\beta -1} \Un_{(0,1)} (x).$ For all $s> -1,$ we have
$$\Delta_{s} = \frac{\beta}{s(\beta + 1)}\lpa 1 - \frac{\Gamma(1/\beta +2)\Gamma(s+2)}{\Gamma(1/\beta + s +2)}\rpa\;\;\mbox{if $s\neq 0$}\qquad\mbox{and}\qquad \Delta_0 = \frac{\beta (\psi(1/\beta +2) -\psi(2))}{\beta +1}\cdot$$
When $\beta = 1/n$ the reciprocal of an integer, the concatenation formula $\psi(z+1) = \psi(z) +1/z$  yields the simple expression
$$\Delta_0 \; =\; \frac{1}{n+1} \lpa \frac{1}{2} +\cdots + \frac{1}{n+1}\rpa$$
and we recover $\Delta_0 = 1/4$ for $\b = 1.$ The computation for $\nabla_{\! s}$ is more involved. Starting from the integral formulation
$$\nabla_{\! s} = \int_0^1 \frac{(1- t^{s+1})(1-t^{1/\beta})}{(1 - t)^2}\, dt \,-\, \frac{1}{\beta + 1}\int_0^1 \frac{(1- t^{s+1})(1-t^{1/\beta + 1})}{(1 - t)^2}\, dt,$$
we obtain
$$\nabla_{\! s}\; =\; \frac{\beta \,(s+1)\,(\psi(1/\beta +2 +s) -\psi(s+2))}{\beta +1}$$
after some simplifications. Observe that again $\nabla_{-1} = \Delta_\infty = 0$ and $\Delta_{-1} = \nabla_\infty = 1/(\b+1) = \esp[X]$ since $z\psi(z) \, \to\, -1$ as $z\to 0$ and $u(\psi(u+z) -\psi(u))\, \to\, z$ as $u\to \infty.$ This example can also be used to compute $\Delta_s$ and $\nabla_{\! s}$ for the exponential distribution $\L,$ which is the limit in law of $ \beta(1- \U^{1/\beta})$ as $\beta\to\infty.$ Using \eqref{Aff}, we obtain 
\begin{equation}
\label{Expo}
\Delta_{s}(\L) \, =\,\frac{\psi(s +2) -\psi(2)}{s}\qquad\mbox{and}\qquad \nabla_{\! s}(\L)\, =\,  (s+1) \psi'(s+2)
\end{equation}
for every $s> -1.$ Observe that both expressions give
$$\Delta_0(\L) \; = \; \nabla_{\!0}(\L) \; =\;\psi'(2) \;= \;\frac{\pi^2}{6} -1,$$
which was recently evaluated in Example 3 of \cite{BBL} by other methods.

\begin{Rem} \label{RAB}
{\em With the notation of the end of Section \ref{Duelle}, the computations of the two previous paragraphs show that 
$$\Diamond_{\U^{1/\beta}}(s) \; =\; s\lpa \frac{\Ga(1/\b +s+2) - \Ga(1/\b +1)\Ga(s+2)}{\Ga(1/\b +s+2) - \Ga(1/\b +2)\Ga(s+2)}\rpa$$
for $s\in (-1,0)\cup (0,\infty),$ and
$$\Diamond_{\U^{1/\beta}}(0) \; =\; \frac{1}{(\b +1) (\psi(1/\beta +2) - \psi(2))}\cdot$$
It is not difficult to show that the mapping $\b\mapsto \Diamond_{\U^{1/\beta}}(s)$ is continuous increasing on $(0,\infty)$ for every $s > -1,$ which implies that $\b\mapsto\varrho(\U^{1/\b})$ is continuous decreasing. Moreover, one has $\Diamond_{\U^{1/\beta}}(1) \to 1$ as $\b\to 0$ and $\Diamond_\U(0) = 1,$ which yields 
$$\lacc \varrho(\U^{1/\b}),\, \b\in (0,1)\racc\; = \; (0,1)$$
by continuity. Similarly, one finds
$${\bar \Diamond}_{\U^{1/\beta}}(s) \; =\; (s+1)(\b(s+1) +1) (\psi(1/\beta +s+2) - \psi(s+2))$$
which is a continuous decreasing function in $\b\in (0,\infty)$ for every $s > -1,$ so that $\b\mapsto{\bar \varrho}(\U^{1/\b})$ is continuous increasing. Moreover, since ${\bar \Diamond}_{\U^{1/\beta}}(s)\to\infty$ as $\b\to 0$ for every $s > -1$ and ${\bar \Diamond}_{\U^{1/\beta}}(s)\to 0$ as $s\to -1$ for every $\b > 0$ together with ${\bar \Diamond}_\U(0) = 1,$ one obtains
$$\lacc \varrho (-\U^{1/\b}) = {\bar \varrho}(\U^{1/\b}),\, \b\in (0,1)\racc\; = \; (-1,0)$$ 
again by continuity. Notice however that for $\b\to\infty,$ one has $\varrho(\U^{1/\b}) \to m = -0.365952..>-1$ and ${\bar \varrho}(\U^{1/\b})\to {\bar m} = 0.389592.. < 1.$ See Figure \ref{RABB} above, and also Remark \ref{RB} below.}
\end{Rem}

\begin{figure}
\centering

\includegraphics[scale =0.8]{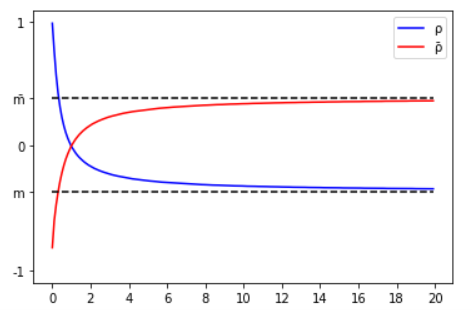}
\caption{Plots of $\b\mapsto \varrho(\U^{1/\b})$ and $\b\mapsto {\bar \varrho}(\U^{1/\b})$ for $\b > 0.$}
\label{RABB}

\end{figure}

\subsection{$X = \U^{-1/\beta} - 1$} 

\label{Lomax}

The density function is $f(x) = \beta (1+x)^{-\beta -1} \Un_{(0,\infty)} (x).$ We recognize the {\em Lomax} distribution, or Pareto distribution of type II, which is the prototype of a power law distribution. We will consider the case with finite expectation only, that is $\beta > 1.$ For $s > -1$ one has the integral formula
$$\Delta_s \;=\; \frac{1}{\beta s} \int_0^1 \frac{x - x^{s+1}}{(1-x)^{1/\beta + 1}} dx \;=\; \frac{1}{\beta}\,\sum_{n\ge 0} \frac{(1/\beta +1)_n}{n! (n+s+2)(n+2)}\cdot$$
The series on the right-hand side is a terminating hypergeometric series ${}_3F_2(1)$ and a consequence of Thomae's relationship and Gauss' formula is then 
$$\Delta_{s}\, = \,\frac{\beta}{s(\beta -1)}\lpa \frac{\Gamma(2-1/\beta)\Gamma(s+2)}{\Gamma(s +2 - 1/\beta)} - 1\rpa\;\;\mbox{if $s\neq 0$}\qquad\mbox{and}\qquad \Delta_0\, =\, \frac{\beta (\psi(2) - \psi(2-1/\beta))}{\beta -1}\cdot$$
This can also be obtained from the computations of the previous paragraph and an analytic continuation at $1/\beta = 0.$ For $\nabla_{\! s}$, the starting point is the integral formula  
$$\nabla_{\! s}\, =\, \int_0^1 \frac{(1- t^{s+1})(t^{-1/\beta} - 1)}{(1 - t)^2}\, dt \,+\, \frac{1}{\beta - 1}\int_0^1 \frac{(1- t^{s+1})(t^{-1/\beta + 1} -1 )}{(1 - t)^2}\, dt$$
for all $s >-1,$ which leads similarly as in Paragraph \ref{Beta2} to
$$\nabla_{\! s}\, =\, \frac{\beta\,(s+1)\,(\psi(s+2) - \psi(s+2 - 1/\beta))}{\beta -1}\cdot$$
One can check that $\nabla_{\! 0} = \Delta_0, \, \nabla_{-1} = \Delta_\infty = 0$ and $\Delta_{-1} = \nabla_\infty = 1/(\b-1) = \esp[X].$ Observe finally from \eqref{Aff} that since $\beta (\U^{-1/\beta} - 1)$ also converges in law to $\L$ as $\beta\to\infty,$ one can also deduce the two formulas in \eqref{Expo} from the Lomax case. 
 
\subsection{$X = 1 - \U^{-1/\beta}$} 

\label{Antilomax}

The density function is $f(x) = \beta (1-x)^{-\beta -1} \Un_{(-\infty, 0)} (x).$ This random variable can be viewed as a negative Lomax. Again, we consider only the case with finite expectation, that is $\beta > 1.$ Observe also that $X_-\in\cL_p$ for every $p <\b$ and that $X\not\in\cL_\b.$ For every $s > -1,$ computations analogous to Paragraphs \ref{Beta1} and \ref{Lomax} give
$$\Delta_s\; =\; \frac{\beta}{(\beta - 1)(\beta (1+s) -1)}\;\;\mbox{if $s> 1/\b -1$}\qquad\mbox{and}\qquad \nabla_{\! s} \; =\; \frac{\beta}{\beta - 1}\lpa \frac{\Gamma(1-1/\beta)\Gamma(s+2)}{\Gamma(s +2 - 1/\beta)} - 1\rpa.$$
One can check that $\nabla_{-1} = \Delta_\infty = 0$ and $\Delta_{1/\b -1} = \Delta_{-1} = \nabla_\infty = \infty = \esp[X] - \min X,$ in accordance with Proposition \ref{Lp}. Observe also that
$$\Delta_0 \;=\; \nabla_{\! 0} \; =\; \frac{\beta}{(\beta -1)^2}$$
and that letting $\beta\to\infty,$ we obtain from \eqref{Aff} and $\beta (1-\U^{-1/\beta})\claw -\L$ the  following formulas for the negative exponential: 
$$\Delta_s(-\L) = \frac{1}{s+1}\qquad\mbox{and}\qquad \nabla_{\! s}(-\L) =\psi(s+2) +\gamma$$
for every $s>-1.$ See also Example 1 in \cite{BBL} for another proof of $\Delta_0(-\L) = 1.$

\begin{Rem} \label{RB}
{\em Similarly as in Remark \ref{RAB} one has
$$\Diamond_{-\U^{-1/\beta}}(s) \; =\; s\lpa \frac{\Ga(s+2-1/\b) - \Ga(1-1/\b)\Ga(s+2)}{\Ga(s+2-1/\b) - \Ga(2-1/\b)\Ga(s+2)}\rpa$$
for $s\in (-1,0)\cup (0,\infty),$ and
$$\Diamond_{-\U^{-1/\beta}}(0) \; =\; \frac{1}{(\b -1) (\psi(2) - \psi(2-1/\b))}$$
which are both continuous decreasing in $\b\in(1,\infty)$ for every $s > -1.$ This shows that $\b\mapsto\varrho(-\U^{-1/\b})$ is continuous increasing on $(1,\infty),$ with $\varrho(-\U^{-1/\b})\to -1$ as $\b \to 1$ since $\Diamond_{-\U^{-1/\beta}}(s)\to\infty$ as $\b\to 1$ for every $s > -1$ and $\Diamond_{-\U^{-1/\beta}}(s)\to 0 $ as $s\to -1$ for every $\b > 1.$ To evaluate the limit as $\b\to\infty$ we observe that $\varrho(-\U^{-1/\b}) = \varrho(\b(1-\U^{-1/\b})) \to\varrho(-\L) = -0.365952..$ which is the unique solution to
$$s\lpa\frac{\psi(s+2) - \psi(1)}{\psi(s+2) -\psi(2)}\rpa\; =\; 1$$
on $(-1,\infty).$ Observe that $\varrho(-\L) = m$ is also the limit of $\varrho(\U^{1/\b}) = \varrho(\b(\U^{1/\b} -1))$ as $\b\to\infty.$ Analogously, we compute
$${\bar \Diamond}_{-\U^{-1/\beta}}(s) \; =\; (s+1)(\b(s+1) -1) (\psi(s+2) - \psi(s+2-1/\b))$$
which is a continuous increasing function in $\b\in (1,\infty)$ for every $s > -1,$ so that $\b\mapsto{\bar \varrho}(-\U^{-1/\b})$ is continuous decreasing on $(1,\infty)$ from 1 to $\varrho(\L) = {\bar m}$ which is the unique solution of 
$$(s+1)^2\psi'(s+2)\; =\;1$$ 
on $(-1,\infty).$ See Figure \ref{RBB} below.}
\end{Rem}

\begin{figure}
\centering

\includegraphics[scale =0.8]{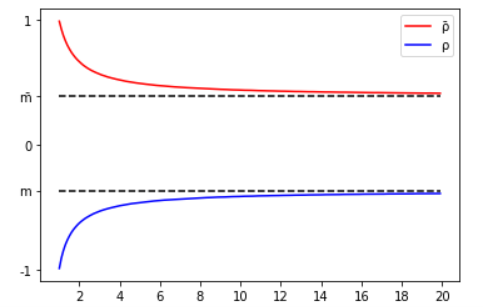}
\caption{Plots of $\b\mapsto \varrho(-\U^{-1/\b})$ and $\b\mapsto {\bar \varrho}(-\U^{-1/\b})$ for $\b > 1.$}
\label{RBB}

\end{figure}
\subsection{$X = \L^{-1/\beta}$} 

\label{Fresh}

The density function is $f(x) = \beta x^{-\beta - 1} e^{- x^{-\beta}} \Un_{\{0,\infty)} (x).$ This is a Fr\'echet distribution, or type II extreme value distribution, and another example of a power law distribution. We consider the case of finite expectation only, that is $\beta > 1.$ For every $s > -1,$ one has 
$$\Delta_s\; = \; \Gamma (1 -1/\beta) \lpa \frac{(s+1)^{1/\beta} - 1}{s}\rpa\;\;\mbox{if $s\neq 0$}\qquad\mbox{and}\qquad \Delta_0\; =\; \frac{\Gamma(1-1/\beta)}{\b}\cdot$$
The formula for $\Delta_0$ was already computed in Example 2.3 of \cite{TP}. For $\nabla_{\! s}$ we get
$$\nabla_{\! s}(X) = \frac{(s+1) \,\Gamma (1-1/\beta)}{\beta}\lpa 1\; +\;\b\,\sum_{n\ge 1} \frac{(-s)_n}{(n+1)!} \lpa \frac{(n+1)^{1/\beta} - 1}{n}\rpa\rpa,$$
which unfortunately does not seem to have a more explicit expression. One has $\nabla_{\!-1} = \Delta_\infty = 0$ and $\Delta_{-1} = \Ga(1-1/\b) =\esp[X]$ but it is not clear from the latter expression that $\nabla_{\! \infty}= \Ga(1-1/\b)$ as well. Considering $\beta (\L^{-1/\beta} - 1)$ with $\beta\to\infty,$ we can compute the $\Delta_s$ and $\nabla_s$ of the Gumbel distribution $-\log \L$, or type I extreme value distribution: 
$$\Delta_{s}( -\log\L) = \frac{\log(s+1)}{s}\qquad\mbox{and}\qquad
\nabla_{\! s}(-\log\L) = (s+1)\lpa 1 + \sum_{n\ge 1} \frac{(-s)_n}{(n+1)!} \frac{\log (1+n)}{n}\rpa$$
for all $s > -1,$ with $\Delta_0(-\log\L) = \nabla_{\! 0}(-\log\L) = 1.$  

\subsection{$X = - \L^{1/\beta}$} 

\label{Weib}

The density function is $f(x) = \beta \vert x\vert^{\beta - 1} e^{- \vert x\vert^{\beta}} \Un_{(-\infty,0)} (x).$ This is a reverse Weibull distribution, or type III extreme value distribution. Similarly as above, one finds 
$$\Delta_s\; = \; \Gamma (1 +1/\beta) \lpa \frac{1- (s+1)^{-1/\beta}}{s}\rpa\;\;\mbox{if $s\neq 0$}\qquad\mbox{and}\qquad \Delta_0\; =\; \frac{\Gamma(1+1/\beta)}{\b}$$
for every $\b > 0$ and $s > -1,$ together with 
$$\nabla_{\! s}(X) = \frac{(s+1) \,\Gamma (1+1/\beta)}{\beta}\lpa 1\; +\;\b\,\sum_{n\ge 1} \frac{(-s)_n}{(n+1)!} \lpa \frac{1-(n+1)^{-1/\beta}}{n}\rpa\rpa.$$
With these expressions, we also retrieve the above $\Delta_s(-\log\L)$ and $\nabla_{\! s}(-\log\L)$ for the Gumbel distribution.

\subsection{$X = \log(\U^{-1} - 1)$} The density function is $e^x/(1+e^x)^2$ over $\rl$ and this random variable is known in the literature as having a logistic distribution. For every $s> -1,$ one finds 
$$\Delta_{s} =\frac{\psi(s +1) +\gamma}{s}\qquad\mbox{and}\qquad \nabla_{\! s} =  \gamma + \psi(s+1) + (s+1) \psi'(s+1),$$
both quantities converging as $s\to 0$ against
$$\Delta_0 \; = \; \nabla_{\! 0} \; =\; \psi'(1)\; =\; \frac{\pi^2}{6}\cdot$$
In the next section, the logistic random variable will appear as the unique maximizer of the cumulative entropy among symmetric random variables in $\cL_2.$

\section{Properties of the range}

In this section, we investigate the closed range of the mapping
$$X\; \mapsto \;\frac{1}{s} \int_\rl F_X(x) (1- F_X(x)^s)\, dx$$
for a fixed $s > -1$ on various subsets, extending some results previously obtained in \cite{BBL, DCL} in the case $s= 0.$ We will deal with the positive case, the case with finite variance, and the symmetric case with finite variance. With an abuse of notation, we will set $\Delta_s(X)$ for the above functional also when $X$ has atoms. Using the standard notation $\sigma_X^{} = \sqrt{{\rm Var} X}$ for $X\in\cL_2$ and introducing the subspaces $\cL_{1,+} = \cL_1 \cap\{ X > 0\}$ and $\cL_{2,{\rm sym}} = \cL_2\cap\{X\elaw -X\},$ we wish to describe the closed sets  
$$\cR_{s,+}\, = \,\overline{\lacc \Delta_s(X)/\esp[X], \, X\in\cL_{1,+}\racc}, \qquad \cR_{s,2}\, =\, \overline{\{ \Delta_s(X)/\sigma_X^{}, \; X\in\cL_2\}}$$
and
$$\cR_{s,{\rm sym}} \, = \, \overline{\{ \Delta_s(X)/\sigma_X^{}, \; X\in\cL_{2,{\rm sym}} \}}$$ 
for $s > -1$ in the positive case resp. for $s > -1/2$ - recall Proposition \ref{Lp} - in the case with finite variance. The above linear normalizations by $\esp[X]$ resp. $\sigma_X^{}$ are necessary in view of the affine relationship \eqref{Aff}. The positive case, which can be viewed as a generalization of the inequality (21) in \cite{DCL} for $s=0,$ is particularly simple.

\begin{Propo}
\label{Posa}
For every $s > -1,$ one has $\cR_{s,+}  = [0,1].$
\end{Propo}

\proof It follows from Proposition \ref{Fall} that $0 \le \Delta_s(X) \le \esp[X]$ for $X\in\DD_+$ and $s > -1.$ The inequalities extend then to all $X\in\cL_{1,+}$ by approximation, since when $X$ has atoms one can construct a sequence $\{X_n, \, n\ge 1\}$ in $\DD_+$ with $X_n\claw X$ viz. $F_{X_n} (x) \to F_X(x)$ at each continuity point $x$ of $F_X,$ which implies 
$$\esp[X_n]\, \to\, \esp[X]\qquad\mbox{and}\qquad \Delta_s(X_n) \,\to\, \Delta_s(X).$$
Considering next $X= \U^{1/\beta}$ for some $\beta > 0,$ we have seen in Paragraph \ref{Beta1} above that 
$$\Delta_s(X)\, =\, \frac{\beta}{(\beta+1)(\beta(1+s)+1)}\qquad\mbox{and}\qquad \esp[X] \; =\; \frac{\beta}{\beta+1}$$
for all $s > -1.$ This implies $\Delta_s(X)/\esp[X] = 1/(\beta (1+s)+1),$ whose closed range is $[0,1]$ as $\beta$ varies from $0$ to $\infty.$ 
\endproof

\begin{Rem}{\em For every $s > -1$ we also have $\overline{\lacc \nabla_{\! s}(X)/\esp[X], \, X\in\cL_{1,+}\racc} = [0,1].$ Indeed, for every $X\in\cL_{1,+}$ one has $0 \le \nabla_{\!s}(X) \le \esp[X]$ with the same argument as above, and we have seen in Paragraph \ref{Beta1} that
$$\frac{\nabla_{\! s}(\U^{1/\b})}{\esp[\U^{1/\beta}]}\; =\; 1 \, -\, \frac{\Gamma(1/\beta +1)\Gamma(s+2)}{\Gamma(1/\beta + s +2)}$$
for $\b > 0,$ whose close range is also $[0,1].$}
\end{Rem}

We next consider the case with finite variance. The following computation improves on all the results of Section 3 in \cite{BBL} - see in particular Theorem 1 therein.

\begin{Propo}
\label{CumRan}
For every $s > -1/2,$ one has $\cR_{s,2} = [0, 1/\sqrt{2s +1}].$   
\end{Propo}

\proof

We first consider the case $s \neq 0,$  starting with the alternative representation
\begin{equation}
\label{Alts}
\Delta_s(X) \; =\; \frac{1}{s}\int_\rl x \lpa (s+1) F^s(x) - 1 \rpa \, dF(x)
\end{equation}
which is valid for all $X\in\DD\cap\cL_2,$ integrating by parts with $ \lcr xF(x) (1- F^s(x))\rcr_{\pm \infty}= 0$ since the square integrability of $X$ implies $x F^{s+1} (x) \le x \sqrt{F(x)} =\sqrt{x^2 F(x)} \to 0$ as $x\to-\infty.$ By the Cauchy-Schwarz inequality, we obtain the required upper bound:
$$\Delta_s^2(X) \, \le \, \frac{{\rm Var} X}{s^2}\, \int_\rl \lpa (s+1) F^s(x) -1\rpa^2  dF(x)\, = \,\frac{{\rm Var} X}{s^2}\, \int_0^1 \lpa (s+1) u^s -1\rpa^2 \, du\, =\, \frac{{\rm Var} X}{2s+1}\cdot$$
Again, this inequality extends to all $X\in\cL_2$ by approximation. Supposing first $s>0$ and setting $X = \U^{1/\beta},$ the computations in Paragraph \ref{Beta1} imply 
$$\frac{\Delta_s(X)}{\sigma^{}_X}\; =\; \frac{\sqrt{\beta(\beta+2)}}{\beta(s+1)+1}$$
which is a unimodal function in $\beta\in (0,\infty)$ from $0$ to $1/(s+1),$ reaching its maximum $1/\sqrt{2s+1}$ at $\beta =1/s.$ By continuity, we get $\cR_{s,2} = [0, 1/\sqrt{2s +1}]$ as required. Supposing next $s\in (-1/2,0)$ and setting $X = 1-\U^{-1/\beta}$ with $\b > 2,$ the computations in Paragraph \ref{Antilomax} imply 
$$\frac{\Delta_s(X)}{\sigma^{}_X}\; =\; \frac{\sqrt{\beta(\beta-2)}}{\beta(s+1)-1}$$
which is a unimodal function in $\beta\in (2,\infty)$ from $0$ to $1/(s+1),$ reaching its maximum $1/\sqrt{2s+1}$ at $\beta =-1/s$ and we can conclude as above. The case $s = 0$ is analogous, starting from the formula
$$\Delta_0(X)\; =\; \int_\rl x \lpa 1 + \log F(x)\rpa\, dF(x)$$
which is a direct consequence of \eqref{Alts} as $s\to 0$ - see also Proposition 2 in \cite{Ra}, and gives the upper bound $\Delta_s(X)/\sigma_X^{} \le 1.$ The open range $(0,1)$ is again described by $\U^{1/\b}$ for $\b\in (0,\infty)$ and we deduce from the end of Paragraph \ref{Antilomax} that the maximum is attained by $-\L$ with $\sigma = \Delta_0 = 1.$   

\endproof

\begin{Rem}{\em (a) In the case $s =1,$ the previous result combined with Remark \ref{Positano} yields 
$$\esp [\vert X - {\tilde X}\vert]\, \le \, \frac{2}{\sqrt{3}}\,\sigma^{}_X$$
for every $X\in\cL_2,$ where ${\tilde X}$ is an independent copy of $X.$ This bound, which we could not locate in the literature, is probably well-known. Observe that applying the Cauchy-Schwarz directly to 
the LHS only leads to $\esp [\vert X - {\tilde X}\vert ] \le  \sqrt{2}\,\sigma^{}_X.$

\medskip

(b) The above proof shows that the maximum of $\cR_{s, 2}$ is attained for every $s > -1/2.$ On the other hand, the maximum of $\cR_{s, +}$ is never attained since Proposition \ref{Fall} implies $\Delta_s(X) < \esp[X]$ for every $X\in \DD_+.$ Observe however that the addition of some constraints may lead to attained upper bounds for cumulative entropies of positive random variables - see Theorem 2 in \cite{Ra}.}
\end{Rem} 

The following result adds a stone to the manifold characterizations of the exponential distribution. This stone is here expressed in terms of the CRE, and seems unnoticed. See \cite{LR} for a related characterization in terms of the relevation transform. 

\begin{Coro} Up to translation, the random variable $\L$ is the unique maximizer of the cumulative residual entropy ${\bar \Delta_0}$ among random variables in $\cL_2$  with unit variance.
\end{Coro}

\proof  By the proof of Proposition \ref{CumRan} and the case of equality in the Cauchy-Schwartz inequality, a random variable $X\in\cL_2$ reaches the maximum of $\cR_{0, 2}$ if and only if there exists some constant $\lambda > 0$ such that $\lambda x  =  (1 + \log F(x))$ for almost every $x\in\rl.$ By continuity, this amounts to $F(x) = e^{\lambda x -1}\Un_{(-\infty, 1/\lambda]}(x)$ that is 
$$X\; \elaw\; \frac{1-\L}{\lambda}\cdot$$ 
Assuming unit variance, we obtain $\lbd = 1,$ which completes the proof since $\Delta_0(X) = {\bar \Delta_0} (-X).$

\endproof

\begin{Rem}{\em It is worth recalling that for the classical Shannon differential entropy 
$$I(X)\; =\; -\int_\rl f_X(x) \log f_X(x)\, dx,$$
the exponential random variable $\L$ maximizes $I(X)$ among absolutely continuous positive distributions with unit expectation, whereas the standard normal random variable $\N$ maximizes $I(X)$ among absolutely continuous real distributions with unit variance, and that 
$$I(\L)\; = \; 1\qquad\mbox{and}\qquad I(\N) \; =\; \frac{1}{2}\, +\, \log\sqrt{2\pi} \, =\, 1.4189..$$
}
\end{Rem}

\medskip

We finally consider the symmetric case with finite variance, which exhibits the most interesting optimal upper bound. In the case $s \in(-1/2, 0)\cup(0,\infty),$ this bound shares some similarities with $\Delta_s (1- \U^s)$ as computed in Paragraph \ref{Beta2}, and we will see during the proof that it is actually reached for $(1-\U)^s - \U^s.$ Observe also the similarity between this bound and the general term of the series appearing in Theorem 5 of \cite{BBL}. Notice finally that the case $s= 0$ improves on Theorems 4 and 5 in \cite{BBL}.

\begin{theo}
\label{sRan}
One has $\cR_{0, {\rm sym}} = [0,\pi/2\sqrt{3}]$ and, for every $s\in (-1/2,0)\cup (0,\infty),$   
$$\cR_{s, {\rm sym}} \; =\; \lcr 0, \frac{s+1}{\sqrt{2 s^2 (2s+1)}}\,\sqrt{{\displaystyle 1\, -\, \frac{\Gamma^2 (s+1)}{\Gamma (2s+1)}}}\rcr.$$
\end{theo}

\proof
We start with the case $s > 0$, using the formula
\begin{equation}
\label{Syms}
\Delta_s(X) \; =\; \frac{s+1}{s}\,\int_0^\infty x \lpa F^s(x) - {\bar F}^s(x) \rpa \, dF(x)
\end{equation}
for all $X\in\DD\cap\cL_{2, {\rm sym}},$ which is a direct consequence of \eqref{Alts}. The Cauchy-Schwarz inequality implies 
\begin{eqnarray*}
\frac{\Delta_s(X)^2}{{\rm Var} X} \; \le \;\frac{(s+1)^2}{s^2}\,\int_0^\infty \lpa F^s(x) - {\bar F}^s(x) \rpa^2 \, dF(x) & = & \frac{(s+1)^2}{2 s^2}\int_{1/2}^1 (x^s - (1-x)^s)^2 \, dx \\
& = & \frac{(s+1)^2}{4 s^2}\int_0^1 (x^s - (1-x)^s)^2 \, dx\\
& = & \frac{(s+1)^2}{2s^2} \lpa \frac{1}{2s+1} \, -\, \int_0^1 x^s(1-x)^s\, dx\rpa\\
& = & \frac{(s+1)^2}{2s^2(2s+1)} \lpa 1 \, -\, \frac{\Gamma^2(s+1)}{\Gamma(2s+1)}\rpa, 
\end{eqnarray*}
as required for the upper bound, which remains valid on the whole $\cL_{2,{\rm sym}}$ by approximation. Observe in passing that the function
$$s\,\mapsto\, \frac{\Gamma^2(s+1)}{\Gamma(2s+1)}$$
has logarithmic derivative $2(\psi(s+1) - \psi(2s+1))$ and is hence unimodal on $(-1/2,\infty)$ with maximum value $1$ at $s=0,$ confirming the positivity of the above RHS. To show the full range, we first consider the case $s > 0$ and introduce the function
\begin{equation}
\label{full}
\phi_s(x) = x^s - (1-x)^s,
\end{equation}
which defines an increasing bijection from $[1/2,1]$ onto $[0,1],$ and for every $\beta \in (0,1]$ the symmetric random variable $X_{s,\beta}$ on $[-1,1]$ with distribution function
$$F_{X_{s,\beta}}(x) = \phi_s^{-1}(x^\beta),\qquad x\in [0,1].$$ From the easily established identity $X_{s,\beta} \elaw \varepsilon \vert X_{s,1}\vert^{1/\beta}$ with $\varepsilon$ an independent random variable such that $\pb[\varepsilon = 1] = \pb[\varepsilon = -1] = 1/2,$ we have
$${\rm Var} X_{s,\beta} \; = \; \esp [ \vert X_{s,1}\vert^{2/\beta}]\; =\; 2\int_0^1 x^{2/\beta} f_s(x)\, dx\; =\; \beta \int_0^\infty e^{-u(1+\beta/2)} f_s(e^{-\beta u/2})\, du$$
where $f_s$ stands for the density of $X_{s,1}.$ As $\beta\to 0,$ this gives the asymptotics
$${\rm Var} X_{s,\beta} \; \sim \; \beta f_s(1) \; =\; \frac{\beta}{s}$$
where the equality is an easy consequence of \eqref{full}. Moreover, it follows from \eqref{Syms} and \eqref{full} that
$$\Delta_s(X_{s,\beta}) \; =\;  \frac{s+1}{s}\,\int_0^1 x^{\beta +1} f_{X_{s,\beta}}(x)\, dx\; =\; \frac{\beta(s+1)}{s}\,\int_0^1 x^{2\beta } f_s(x)\, dx\;\sim\; \frac{\beta(s+1)}{2s}$$
as $\beta\to 0,$ which implies
$$\frac{\Delta_s(X_{s,\beta})^2}{{\rm Var} X_{s,\beta}}\, \sim\,  \frac{\beta (s+1)^2}{4 s}\, \to \, 0\qquad\mbox{as $\beta\to 0.$ }$$
Finally, for $\beta = 1$ the above computation also gives
\begin{eqnarray*}
\frac{\Delta_s(X_{s,1})^2}{{\rm Var} X_{s,1}} \; = \;\frac{(s+1)^2 {\rm Var} X_{s,1}}{4 s^2} & = &\frac{(s+1)^2}{2 s^2}\, \int_0^1 x^2 f_s(x) dx\\
& = & \frac{(s+1)^2}{2 s^2}\int_{1/2}^1 (x^s - (1-x)^s)^2 \, dx \\
& = & \frac{(s+1)^2}{2s^2} \lpa \frac{1}{2s+1} \, -\, \frac{\Gamma^2(s+1)}{\Gamma(2s+2)}\rpa,
\end{eqnarray*}
showing that the upper bound is attained by $(1-\U)^s - \U^s$, and we can conclude by continuity. The case $s \in (-1/2,0)$ is analogous, using the function $\psi_s(x) = (1-x)^s - x^s$ which defines an increasing bijection from $[1/2,1]$ onto $\rl^+,$ and the same maximizing random variable $(1-\U)^s - \U^s$ which is here unbounded. We omit details. Finally for the case $s = 0,$ we use the formula 
$$\Delta_0(X)\; =\; \int_0^\infty\!\! x \,\log\lpa\frac{F(x)}{{\bar F} (x)}\rpa dF(x),$$
which is a direct consequence of \eqref{Syms} as $s\to 0.$ The Cauchy-Schwarz inequality leads here to 
\begin{eqnarray*}
\Delta_0(X)^2 & \le & \lpa \int_0^\infty x^2 \, dF(x)\rpa\times\lpa \int_0^\infty \log^2\lpa\frac{F(x)}{{\bar F} (x)}\rpa dF(x)\rpa\\
& = & \frac{{\rm Var X}}{2}\,\times\lpa \int_{1/2}^1 \log^2(x/(1-x))\,dx\rpa \; =\; \frac{\pi^2 \,{\rm Var X}}{12} 
\end{eqnarray*}
which is the required upper bound. This upper bound is attained by the logistic random variable $X = \log (\U^{-1} -1),$ and the full range is described by its symmetric powers, as above. 

\endproof

\begin{Rem}{\em (a) The maximizing random variable $X_{s,1} =  (1-\U)^s - \U^s$ is bounded for $s > 0$ and unbounded with heavy tails for $s \in (-1/2,0).$ It is remarkable that the same dichotomy occurs for the so-called Tsallis or $q-$Gaussian random variable, which depends on some parameter $q < 3$ and can be constructed as
$$\T_q\; =\; \varepsilon \,\sqrt{\frac{(3-q) (1- \U^{1-q})}{1-q}}$$
for $q\neq 1$ with independent $\varepsilon$ and $\U$ as in the preceding proof, and as the Gaussian limit $\T_1 =\N.$  Indeed this symmetric random variable $\T_q,$ which is a maximizer of the Tsallis entropy
$$I_q(X) \; =\; \frac{1}{q-1} \lpa 1\, -\, \int_\rl f_X(x)^q \, dx\rpa$$
under appropriate constraints - see \cite{PT} and especially formula (9) therein, is bounded for $q < 1$ and unbounded with heavy tails for $1 < q < 3.$ In this respect, the random variable $X_{s,1}$ may be called the $s-$Logistic random variable.
 
\medskip

(b) The density of the $s-$Logistic random variable $X_{s,1}$ does not seem to have an explicit character in general, except for $s=1,2$ where $X_{1,1}$ and $X_{2,1}$ are uniform on $[-1,1],$ and for $s=1/2$ where $X_{1/2,1}$ has density
$$\frac{1-x^2}{\sqrt{2-x^2}}\, \Un_{(-1,1)} (x).$$
This contrasts with the explicit density of the above random variables $\T_q,$ which reads $e^{-x^2/2}/\sqrt{2\pi}$ for $q = 1$ and
$$C_q \lpa 1 + \lpa\frac{q-1}{3-q}\rpa x^2\rpa_+^{\frac{1}{1-q}}$$
for $q\neq 1,$ where $C_q$ is the normalizing constant which can be computed in closed form.

\medskip

(c) The non-increasing character of $s\mapsto\Delta_s(X)$ implies that $\{\cR_{s, {\rm sym}}, \, s > -1/2\}$ is a non-increasing family of intervals, expanding to $\rl^+$ as $s\to -1/2$ and shrinking to $\{0\}$ as $s\to\infty.$ Considering the upper bound shows the non-trivial fact that the mapping
$$s\;\mapsto\; \frac{(s+1)^2}{2s^2(2s+1)}\lpa 1\, -\, \frac{\Gamma^2(s+1)}{\Gamma(2s+1)}\rpa$$
decreases on $(-1/2,\infty)$ from $\infty$ to 0. One might ask if this mapping is not completely monotone. See \cite{AB} for several completely monotonic functions related to the Gamma function. }
\end{Rem}

Let us now mention the following characterization of the logistic distribution as a maximizer of the cumulative entropy, a noteworthy counterpart to that of the exponential distribution as a maximizer of the cumulative residual entropy. 

\begin{Coro} Up to translation, the rescaled logistic random variable $\frac{\sqrt{3}}{\pi}\,\log(\U^{-1} -1)$ is the unique maximizer of the cumulative entropy among random variables  in $\cL_{2,{\rm sym}}$ with unit variance.
\end{Coro}

\proof Similarly as above, the proof of Theorem \ref{sRan} and the case of equality in the Cauchy-Schwarz inequality show that a random variable $X\in\cL_{2, {\rm sym}}$ reaches the maximum of $\cR_{0, {\rm sym}}$ if and only if there exists some constant $\lambda > 0$ such that
$$\lambda\, x \; = \; \log\lpa\frac{F_X(x)}{{\bar F}_X (x)}\rpa $$
for almost every $x\in\rl^+,$ which by symmetry and continuity amounts to   
$$\frac{F_X(x)}{1 - F_X (x)}\; =\; e^{\lambda x}$$
for every $x\in\rl,$ that is $X \elaw \frac{1}{\lambda} \log(\U^{-1} -1).$ Finally, the constraint of having unit variance gives $\lbd = \pi/\sqrt{3}$ by the known formula
$${\rm Var} (\log(\U^{-1} -1))\; = \; \frac{\pi^2}{3}\cdot$$
\endproof

\begin{Rem}{\em The cumulative entropy of the standard Gaussian random variable, which maximizes the Shannon entropy among real, symmetric or not, absolutely continuous random variables with unit variance, can be rewritten as
$$\frac{1}{\sqrt{2\pi}} \int_0^\infty e^{-u} \,\log \lpa 2e^{u+\varphi(u)} -1\rpa\,du$$
where $\varphi(u) = -\log\esp[e^{-uA_{1/2}^{-1}}]$ and $A_{1/2}$ is the arcsine random variable of Remark \ref{Senans} (a). However, this integral does not seem to be computable in closed form. Some simulations give an approximate value $\Delta_0(\N) = 0.9033...$ which is smaller than but close to $\pi/2\sqrt{3} = 0.9068...$} 
\end{Rem}

\medskip

We conclude this paper with a non-trivial inequality for the classical Gamma function, which is in the case $s\in(-1/2,0) \cup (0,1)$ a consequence of Theorem \ref{sRan}. 

\begin{Coro}
\label{GIneq}
For every $s\ge -3/2,$ one has
\begin{equation}
\label{RueGama}
\frac{\Gamma^2(s+2)}{\Gamma(2s +1)} \; \ge \; 1 + 2s - s^2
\end{equation}
Moreover, the inequality is strict except at $s =0,1.$
\end{Coro}

\proof
The equality is plain for $s=0$ or $s=1,$ and the strict inequality is also straightforward for $s \in [-3/2,-1]\cup [1+\sqrt{2},\infty)$ since then the RHS is non-positive and the LHS is non-negative. The strict inequality for $s\in (-1,-1/2]$ is obtained directly from the equivalent formulation 
\begin{equation}
\label{Tumulte}
\frac{\sqrt{\pi}(s+1)\Gamma(s+2)}{4^s \Gamma(s+1/2)}\; \ge\; 1 + 2s - s^2
\end{equation}
which is given by the Legendre duplication formula, and holds true for $s=-1/2$ since the left-hand side is zero. If $s\in (-1,-1/2)$ the inequality \eqref{Tumulte} amounts to
$$ \lpa\frac{\sqrt{\pi}\,4^t\, \Gamma(3/2-t)}{\Gamma(1-t)}\rpa t(1-2t)\; \le\; 1/4 + 3t + t^2$$
for $t =-1/2-s\in(0,1/2),$ whose LHS is bounded from above by $\pi t(1-2t)$ by log-convexity of the Gamma function, and a trinomial analysis shows that $\pi t(1-2t) < 1/4 + 3t + t^2$ for all $t\in\rl.$ 

We next consider the strict inequality for $s \in (1, 1+\sqrt{2}]$. Taking the logarithmic derivatives on both sides, we are reduced to show that
$$\psi(s+2) \, -\, \psi(2s+1)\; > \; \frac{1-s}{1+2s -s^2}$$
for all $s\in (1,1/(\sqrt{2} -1)).$ This amounts to
$$\frac{1}{1+2s-s^2}\, -\, \sum_{n\ge 1} \frac{1}{(n+2s)(n+s+1)}\; > \; 0,$$
which holds true since the LHS equals $9/4 - \pi^2/6 > 0$ at $s =1$ and increases on $(1,1+\sqrt{2}].$ 

We finally show the strict inequality for $s\in (-1/2,0)\cup (0,1),$ which cannot seem to be handled neither directly nor with classical monotonicity or convexity arguments. Instead, we consider the equivalent formulation  
$$\frac{s+1}{\sqrt{2 s^2 (2s+1)}}\,\sqrt{{\displaystyle 1\, -\, \frac{\Gamma^2 (s+1)}{\Gamma (2s+1)}}}\; < \; \frac{1}{\sqrt{2s +1}}$$
for $s\in (-1/2,0)\cup(0,1).$ By the proof of Theorem \ref{sRan} and using the notation $\Delta_s^\sigma (X) = \Delta_s(X)/\sigma_X^{}$ for $X\in\DD_2,$ this is tantamount to $\Delta_s^\sigma ((1-\U)^s - \U^s) <  \Delta_s^\sigma (\U^s).$ But it is clear by definition that 
$$\Delta_s^\sigma ((1-\U)^s - \U^s)\; \le \; \max\{\Delta_s^\sigma (X), \, X\in \DD\cap\cL_2\}\; =\; \Delta_s^\sigma(\U^s)$$
and by uniqueness that the inequality is strict, since any affine transformation of $\U^s$ is not symmetric for $s\in (0,1)$ and hence cannot be distributed as $(1-\U)^s - \U^s.$ 
 
\endproof

\begin{Rem} {\em (a) The inequality \eqref{Tumulte} amounts to 
$$\frac{\Gamma(s+1)}{\Gamma(s+1/2)}\; \ge\; \frac{4^s (1+2s- s^2)}{\sqrt{\pi} (1+s)^2}$$ 
for all $s > -1/2,$ whose right-hand side can be shown to be greater than $\sqrt{s+1/4}$ for all $s\in [-1/4,1]$ by an elementary monotonicity argument. In particular, \eqref{RueGama} can be viewed as an improvement on Watson's inequality. We refer to \cite{QL} and the references therein for a collection of classical inequalities for the Gamma function including Watson's, none of which seems to imply \eqref{RueGama} directly. We observe that \eqref{RueGama} is rather sharp on $(0,1)$: simulations show that 
\begin{equation}
\label{Gamm}
s\;\mapsto\; \varphi(s) \, =\, \frac{\Gamma^2(s+2)}{\Gamma(2s +1)} - 1 - 2s + s^2
\end{equation}
is unimodal on $(0,1)$ from 0 to 0, with a small maximum value 0.0172.. attained at $s = 0.4671..$ See Figure \ref{Tux} above.

\begin{figure}
\centering

\includegraphics[scale =0.4]{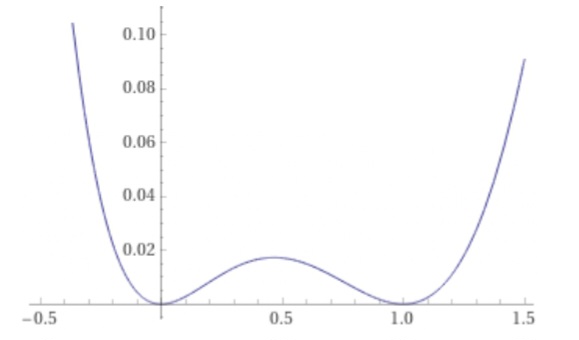}
\caption{Plot of $s\mapsto \varphi(s)$ for $s \in [-1/2,3/2].$}
\label{Tux}

\end{figure}

(b) The inequality \eqref{RueGama} becomes false as $s \to -2$ because the left-hand side tends to $-\infty.$ On the other hand, the derivative of \eqref{Gamm} equals
$$2\lpa \frac{\Gamma^2(s+2) \,\lpa \psi(s+2) - \psi(2s+1)\rpa}{\Gamma(2s +1)}\, -\, s\, -\, 1\rpa\, >\, 0$$
for $s\in(-2,-3/2)$ since $\Ga(2s+1) < 0$ and $\psi(s+2) < \psi(2s+1)$ on this interval, as can be easily checked. Putting everything together shows that \eqref{RueGama} actually holds for all $s \ge s_* = -1.6609...$ which is the unique root of \eqref{Gamm} on $(-2,-3/2).$ See Figure \ref{Tuxi} below.

\begin{figure}
\centering

\includegraphics[scale =0.4]{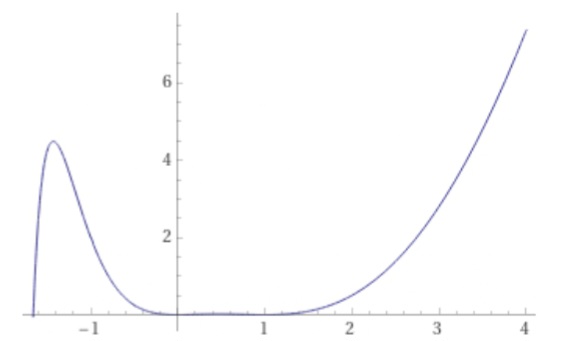}
\caption{Plot of $s\mapsto \varphi(s)$ for $s \in [s_*,4].$}
\label{Tuxi}

\end{figure}

(c) It is interesting to compare \eqref{RueGama} on $(0,1)$ with the inequality 
$$B(x,y) \; > \; \frac{x+y}{xy}\lpa 1\, -\, xy \lpa \frac{2}{x+y+1}\wedge \frac{1}{x+y}\rpa\rpa$$
for all $x,y\in (0,1)$ obtained recently in \cite{ZW} - see also the references therein, where $B(x,y)$ is the classical Beta function. Setting $x = y =s$ therein leads indeed to  
$$\frac{\Gamma^2(s+2)}{\Gamma(2s +1)} \; > \; 1 + 2s - s^2 - \lpa\frac{2s^4}{2s+1}\wedge \frac{s(s-1)^2}{2}\rpa,$$
which is slightly less sharp than \eqref{RueGama}.}

\end{Rem}

\end{document}